\documentclass[12pt]{article}

\usepackage{amssymb}
\usepackage{amsmath}
   
\newtheorem{prop}{Proposition}[section]
\newtheorem{theor}{Theorem}[section]
\newtheorem{lemma}{Lemma}[section]

\newtheorem{remark}{Remark}[section]

\def\ox{\overline x}
\def\oy{\overline y}
\def\eps{\varepsilon}

\def\e{\varepsilon}
\def\11{1\!\!1}

\newcommand{\sn}{~\hfill\hbox{$\square$}}

\def\oy{\overline y}

\def\Ktu{{\tilde K_1(Y,\delta)}}
\def\Ktd{{\tilde K_2(Y,\delta,\eta)}}

\def\Dx1 #1{{{\partial #1}\over{\partial x_}}}

\def\rec#1{{(\ref{#1})}}
\def\e{{\varepsilon}}

\def\r{\right}

\def\BUC{B\kern-.5pt U\kern-1pt C}

\def\preqed{\hbox to 8pt{\vbox{\hsize 6pt\hrule width5pt\noindent\vrule
height7pt \hskip 4pt\vrule height 7pt\hrule width 5pt}\hfil}}

\def\text{\textstyle}

\def\tp{{\cal T}(p)}

\renewcommand{\r}{{\mathbb R}}
\newcommand{\R}{\r^{n}}

\begin{document}

\title{\bf Local $C^{0,\alpha}$ Estimates for Viscosity Solutions of
Neumann-type Boundary
Value Problems}
\author{Guy\ Barles$^{(1)}$ \&
Francesca Da Lio$^{(2)}$}

\addtocounter{footnote}{1} \footnotetext{Laboratoire de
Math\'ematiques et Physique Th\'eorique.  Universit\'e de Tours.
Facult\'e des Sciences et Techniques, Parc de Grandmont, 37200
Tours, France.} \addtocounter{footnote}{1}
\footnotetext{Dipartimento di Matematica.  Universit\`a di Torino.
Via Carlo Alberto 10,  10123 Torino, Italy.}

\date{ }

\maketitle

\begin{abstract}
In this article, we prove the local $C^{0,\alpha}$ regularity and
provide $C^{0,\alpha}$ estimates for viscosity solutions of fully
nonlinear, possibly degenerate, elliptic equations associated to linear
or nonlinear Neumann type boundary conditions. The interest of these
results comes from the fact that they are indeed regularity results
(and not only a priori estimates), from the generality of the equations
and boundary conditions we are able to handle and the possible
degeneracy of the equations we are able to take in account in the case
of linear boundary conditions.
\end{abstract}

\maketitle
   \date{}
\tableofcontents
\section{Introduction}

In this article, we are interested in the local $C^{0,\alpha}$
regularity of viscosity solutions of nonlinear Neumann boundary
value problems of the form
\begin{equation}\label{np} \left\{\begin{array}{ll}
   F(x,u,Du,D^2u)=0 & \mbox{in $O$,} \\
\displaystyle G(x,u,Du)=0& \mbox{on $\partial O $,} \\
    \end{array} \right.
    \end{equation}
where $O\subset\R$ is a smooth domain, $F$ and $G$ are, at least,
real-valued continuous functions defined respectively  on $\overline O
\times \r \times\R\times {\cal{S}}^n$ and $\partial
O\times\r\times\r^n, $ where ${\cal{S}}^n$ denotes the space of real,
$n\times n$, symmetric matrices.    The solution $u$ is a scalar
function and $Du$, $D^2u$ denote respectively its gradient and Hessian
matrix. More precise assumptions on $O,F,G$ are given later on.

We recall that the boundary condition $G=0$ is said to be a nonlinear
Neumann boundary condition if the function $G$ satisfies the following
conditions
\begin{enumerate}
\item[{\bf(G1)}] For all $R>0$, there exists $\mu_{R}>0$ such
that, for every $(x,u,p)\in \partial O\times[-R,R]\times\r^n,$
and $\lambda>0$, we have
\begin{equation}\label{croiss}
G(x,u,p+\lambda n(x))-G(x,u,p)\ge \mu_{R}\lambda\,,
\end{equation}
\hspace{-\leftmargin}where $n(x)$  denotes the unit outward normal
vector to $\partial
O$ at $x \in \partial O$.
\item[{\bf(G2)}] For all $R>0$ there
is  a constant $K_{R} >0$ such that, for all $x,y\in
\partial O,$ $p,q\in\r^n$, $u,v\in [-R,R] $, we have
\begin{equation}
|G(x,u,p)-G(y,v,q)|\le K_R\left[(1+|p|+|q|)|x-y|+|p-q|+|u-v|
\right]\,.
\end{equation}
\end{enumerate}

The main examples of boundary conditions we have in mind are the
following~: first, linear type boundary conditions like oblique
derivative boundary conditions, in which  $G$ is given by
\begin{equation}\label{odbc}
   G(x,u,p)=\langle p,\gamma(x)\rangle+g(x) \,,
\end{equation}
where $\gamma\colon\partial O\to\r^n$ is a bounded, Lipschitz
continuous vector field such that
$$\langle \gamma(x), n(x)\rangle \ge \beta>0~~~\mbox{for all   $x\in
\partial O$,}$$ and $g \in C^{0,\beta}(\partial O)$ for some
$0<\beta<1$. Here and below, ``$\langle p , q\rangle$'' denotes  the
usual scalar product of the vectors $p$ and $q$ of $\R$.

Next nonlinear boundary conditions~: the first example is capillarity
type boundary conditions for which $G$
is given by
\begin{equation}\label{capillarity}
   G(x,u,p)=\langle p, n(x)\rangle-
\theta(x)\sqrt{1+|p|^2} \,,
\end{equation}
where $\theta\colon\partial O\to \R $ is a Lipschitz scalar function,
such that $|\theta(x)|<  1$ for every $x \in \partial O$. A second
example is the
boundary condition arising in the optimal control of processes with
reflection when there is control on the reflection, namely
\begin{equation}\label{crefl}
G(x,u,p)=\sup_{\alpha\in A
}\{\langle \gamma_{\alpha}(x), p\rangle + c_{\alpha}(x) u
-g_{\alpha}(x)\} ,
\end{equation}
where $A$ is a compact metric space, $\gamma_{\alpha}\colon\partial
O\to\r^n$ are Lipschitz
continuous vector fields such that $\langle\gamma_{\alpha}(x),
n(x)\rangle\ge \beta>0$ for all $x \in \partial O$,  and $ c_{\alpha},
g_{\alpha}\colon\partial O\to \r $ are
Lipschitz continuous, scalar functions.

We are going to show that, under suitable assumptions on $O$, $F$, $G$,
any continuous viscosity solution of (\ref{np}) is in
$C^{0,\alpha}_{loc}(\overline O)$ for some $0<\alpha <1$, with an
estimate on the local $C^{0,\alpha}$-norm of $u$.  These results are
indeed regularity results (and not only a priori estimates);  this is
their first main advantage. But they are also valid, in the case of
linear boundary conditions, for possibly degenerate equations, a second
original feature. The counterpart is that the regularity properties we
have to impose on $O$, $F$, and $G$, are stronger than in the case
of a priori estimates where the solution $u$ is assumed to be either in
$C^2( O)\cap C^1(\overline O)$ or at least in $W^{2,n}(O)\cap
C^1(\overline O) $.

The classical a priori estimates for this type of problems are in fact
proved for linear equations and extended to fully nonlinear
equations by a simple linearization procedure described in
Lieberman \& Trudinger\cite{LT}; this linearization requires the
regularity of the solution.  For linear equations, the classical
results are obtained under rather weak assumptions on the
coefficients of the operators in the equation and in the boundary
condition. To the best of our knowledge, the first results in this
direction are the ones of Nadirashvili \cite{n1,n2} for linear,
uniformly elliptic equations with $L^\infty$ coefficients associated to
oblique derivative boundary condition of the form
\begin{equation}\label{oblder}
\langle Du,\gamma(x)\rangle +a(x) u(x)+g(x)=0\quad\mbox{in $\partial
O$}.
\end{equation}
He first proves them with a continuous direction of reflection
$\gamma$ and for Lipschitz domains, and then for a
direction of reflection in $L^\infty$ for $C^2$ domains.  Results
in this direction were also obtained in the 80's by Liberman
\cite{Li87} by different methods. Recent improvements on the regularity
of
the coefficients of the equation (which can be assumed to be
$L^p$ for $p$ large enough or in $L^n$), were obtained by Kenig \&
Nadirashvili\cite{kn} and Lieberman \cite{lib}.

The case of fully nonlinear equations was first considered in Lions
\& Trudinger \cite{LionsT} who show the existence of a smooth
classical solution in $C^2(O)\cap C^{1,1} (\bar O)$ for
Hamilton-Jacobi-Bellman equations with smooth coefficients and
directions of reflection.  As mentioned above, in Lieberman \&
Trudinger\cite{LT}, the case of fully nonlinear equations is
considered in a more systematic way but most of the results are
obtained by a linearization procedure and are based on results for
linear equations; it is worth pointing out that, in general, the
passage from a priori estimates to regularity results
requires the existence of smooth enough solutions for a sequence of
approximate problems (and even uniqueness for the problem itself),
and is often only valid for convex or concave equations.

Our approach is based on classical viscosity solutions (not $L^p$
viscosity solutions): for a detailed presentation of the theory of
viscosity solutions and of the boundary conditions in the viscosity
sense, we refer the reader to the ``Users'guide" of Crandall, Ishii
and Lions \cite{cil} and the book of Fleming and Soner \cite{fs}, while
the books of Bardi and Capuzzo Dolcetta \cite{bcd} and Barles
\cite{bb} provide an introduction to the theory in the case of
first-order equations.

Clearly, this approach requires more regularity properties for the
operator $G$  which has to be locally Lipschitz with respect to
its variables, for the domain $O$ which has to be assumed to be
$C^2$ and for the equation ($F$ has to be continuous). Its
advantage is that, in the case of oblique boundary conditions of
the form \rec{oblder},  we require only $F(x,u,p,M)$ to be
nondegenerate (in a sense precised below) in one direction which
depends, near the boundary, only on $p$ and $\gamma$. Whereas in
the case of more general nonlinear boundary conditions we require
$F$ to be uniformly elliptic. To prove such results, we use
systematically an idea introduced in Ishii \& Lions \cite{il}
which has already been used to obtain interior, local regularity
(or global regularity) in \cite{Ba1} and Barles \& Souganidis
\cite{BS}.

In this paper, for technical reasons, we treat separately the ``linear
case'', i.e. typically the case of oblique
derivative boundary condition where the operator $G$ is linear with
respect to $u$ and $p$ and the ``nonlinear'' case where $G$ is not
linear. A surprising fact in the linear case - and maybe our result is
not optimal in this direction - is  that the
assumptions on $F$, and in particular the ellipticity one, depends on
$\gamma$. We were unable to remove this dependence.

In the case of Neumann boundary condition i.e. when $\gamma(x) =
n(x)$,  our ``strong ellipticity condition'' can be written
formally as
\begin{equation}\label{egd}
  \frac{\partial F}{\partial M}(x,u,p,M) \leq -
\lambda \hat p \otimes \hat p\quad \hbox{for almost every $(x,u,p,M)$},
\end{equation}
where $\lambda >0$ and where, here and below, the notation $\hat
p$ stands for $\displaystyle\frac{p}{|p|}$. This condition is the
natural requirement for the interior $C^{0,\alpha}$ regularity to
hold and it allows to extend the results up to the boundary.

Next, if the direction of reflection $\gamma$ is $C^2$, then a
classical property which is used in Lions \cite{LiNeu} (see also Lions
\& Sznitman \cite{LiSz}) is the existence of a $C^2$ function $A(x)$,
taking values in the set of non-negative symmetric matrices and such
that $A(x) \gamma (x) =  n(x)$ for any $x\in \partial O$. In this case,
we have to require that the above ``strong ellipticity condition'' is
valid but replacing in (\ref{egd}) $\hat p$ by $\widehat{A^{-1}(x) p}$.
Since $A$ is not unique, this assumption is admittedly not completely
satisfactory.

Finally, if $\gamma(x) \neq n(x)$ is just Lipschitz continuous,
then we have again to assume ``strong ellipticity condition''  of
$F$ in the $\widehat{A^{-1}(x) p}$ direction where again $A(x)
\gamma (x) =  n(x)$ but, here, $A(x)$ is just Lipschitz
continuous and this creates technical difficulties..

We mention that, both in the linear and nonlinear
case, we prove the regularity result by assuming that $G$ does not
depend on $u$. Indeed one can always reduce to this case by a
suitable change of variable that we show later on.

The proofs of these results rely on the constructions of suitable
test-functions inspired by the test-functions built for proving
uniqueness results : in the case of Neumann or regular oblique
derivatives problems, the corresponding uniqueness results were
proved by Lions \cite{LiNeu} (see also \cite{cil}) and in the case
of nonlinear Neumann-type boundary boundary condition in
\cite{BaNeu}. It is worth pointing out anyway that the
construction in the case of Lipschitz continuous $\gamma$'s, which
is the difficult case, takes a completely different form here.

It is worth mentioning also the results of Ishii \cite{I4} proved,
in the case of nonlinear Neumann boundary conditions, under weaker
assumptions on $O$ but stronger assumptions on the boundary
condition than in \cite{BaNeu}; our approach requires more
regularity of the boundary and therefore we do no use the
test-function built in \cite{I4}.

This paper is organized as follows~: in Section 2, we state our
regularity results both in the linear (Subsection 2.1) and in the
nonlinear case (Subsection 2.2) and we provide the main proofs.
Such proofs rely on the constructions of a suitable test-functions
which are different in the linear and nonlinear case: these
constructions are given in Section~\ref{CTF}.
  It is worth pointing out anyway that, despite most of the arguments
are common in these two cases,
  the conclusion is a little bit different because of the particular ``
ellipticity conditions" used in these two cases.

\section{The Local $C^{0,\alpha}$ Estimates}

In this section, we state and prove the local $C^{0,\alpha}$ regularity
of the
solutions of the problem (\ref{np}) both in the linear and nonlinear
cases. As pointed out in the Introduction, these two cases requires
slightly different assumptions. We first introduce the assumptions
which are common of both cases. First, for the domain $O$, we require

\begin{enumerate}
\item[{\bf (H1)}]~~{\bf (Regularity of the boundary)} $O$ is a
domain with a $C^2$--boundary.\par
\end{enumerate}
This assumption on $O$ implies the existence of an $\R$-neighborhood
${\cal{V}}$ of $\partial O$ such that the signed distance function
$d$ which is positive in $O$ and negative in $O^c$ is in
$C^2({\cal{V}}).$ We still denote by $d$ a $C^2$-extension of the
signed distance function to $\R$ which agrees with $d$ in ${\cal{V}}$
and we use below the notation $n(x) = - Dd(x)$ even
if $x$ is not on the boundary.

The  ``strong ellipticity" conditions on $F$ are different in the
linear and nonlinear cases but
the following natural growth condition on $F$ is, on the contrary, the
same

\begin{enumerate}
\item[{\bf (H2)}]~~{\bf (Growth Condition on $F$)} For any $R>0$,
there exist positive constants $C_1^R$, $C_2^R$, $C_3^R$ and
functions $\omega_1^R$, $\omega_2^R$, $\varpi^R : \r^+ \to \r$ such that
$\omega_1^R (0+)=0$, $\omega_2^R(r)=O(r)$ as $r\to + \infty,$ $\varpi^R
(t) \to 0 $ as $t\to+\infty$,   and for
any $x,y \in \overline O$, $-R\leq u,v \leq R$, $p,q \in \R$, $M
\in {\cal{S}}^n$ and $K>0$ \begin{eqnarray*} F(x,u,p,M)-F(y,v,q,M
+ KId)     & \leq &  \omega_1^R(|x-y|(1+ |p|+|q|)\\ &+& \varpi^R (|p|
\wedge |q|)|p-q|) ||M|| \\
& & + \omega^R_2(K) +
C_1^R  + C_2^R(|p|^2+|q|^2) \\
       &  & +C_3^R|x-y|(|p|^3+|q|^3)\;  ,
\end{eqnarray*}
where $|p|
\wedge |q| = \min (|p|, |q|)$.
\end{enumerate}

In the sequel, $K$ always denotes a positive constant which may vary
from line to line, depends only on the
data of the problem and is, in particular, independent on the small
parameters we are going to introduce.

\subsection{The case of linear boundary conditions}

In this subsection we examine the case when $G$ is linear with
respect to $p$, namely it is of the form \rec{oblder}.

The main additional assumptions on $F$ and $G$ are the following.
\begin{enumerate}
\item[{\bf (H3a)}] {\bf  Oblique-derivative boundary
condition and ellipticity~:} there exists a \linebreak Lipschitz
continuous function $A\colon\overline O\to  {\cal{S}}^n$
with $A\ge c_0Id,$ for some $c_0>0$ such that $A(x)\gamma(x)=n(x)$ for
every $x\in
\partial O$, and for any $R>0$, there exist $L_R, \lambda_{R}
0$ such that, for all $x\in \overline O,$ $|u|\le R,$
$|p|>L_{R}$ and $M,N\in {\cal{S}}^n$ with $N\geq 0$, we have
\begin{equation}
   \label{ndpsod}
   F(x,u,p,M + N )-F(x,u,p,M) \le -\lambda_{R} \langle N
\widehat{A^{-1}(x)p},
\widehat{A^{-1}(x)p}\rangle + o(1)||N|| \, ,
\end{equation}
where $o(1)$ denotes a function of the real variable $|p|$ which
converges to $0$ as $|p|$ tends to infinity.
\end{enumerate}

Before stating the assumption on the boundary condition, we want
to point out that the existence of such $A$ is really an
assumption in a neighborhood of $\partial O$, then, under suitable
assumptions on $F$, $A$ can be extended to $\overline O$ and even
to  $\R$.

For the boundary condition , we require
\begin{enumerate}
\item[{\bf (H4)}] {\bf (Regularity of the boundary condition)}
The functions $\gamma$ and $a$ in (\ref{oblder}) are Lipschitz
continuous on $\partial O$, $\langle \gamma(x), n(x)\rangle\ge \beta>0$
for any $x\in \partial O$  and  $g$ is
is in $C^{0,\beta}_{loc}(\partial O)$ for some  $0<\beta \leq 1$.
\end{enumerate}

Our result is the
\begin{theor}\label{theorregul}
Assume {\bf (H1)-(H2)-(H3a)-(H4)}.  Then every continuous
viscosity solution $u$ of \rec{np} with $G$ given by
(\ref{oblder}) is in $C^{0,\alpha}_{loc}(\overline O)$ for any
$0<\alpha< 1$ if $\beta = 1$ and  with $\alpha=\beta$ if
$\beta<1$. Moreover the $C^{0,\alpha}_{loc}$--norms of $u$ depend
only on $O$, $F$, $\gamma$, $a$, $g$ through the constants and
functions appearing in {\bf  (H2)-(H3a)}, the local
$C^{0,1}$--norm of $\gamma$ and $a$, the local $C^{0,\beta}$--norm
of $g$ and the local $C^2$--norm of the distance function of the
boundary including the modulus of continuity of $D^2 d$.
\end{theor}

Despite we make a point here to have a rather unified result as we do
it for the proof, Theorem~\ref{theorregul} contains clearly three cases
which are rather different from the technical point of view.

\noindent 1. The homogeneous Neumann boundary condition $\displaystyle
\frac{\partial u}{\partial n} = 0$ or more generally
$\displaystyle \frac{\partial u}{\partial n} + g(x)=0$ when $g$ is
a $C^2$-function, is the simplest case. Of course, one can take $A$ as
being the identity matrix and the construction of the test-function
does not require the heavy regularization procedures we use in
Section~\ref{CTF}.

It is worth pointing out
here that, in the construction of the test-function, the two terms
$\displaystyle\frac{\partial u}{\partial n}$ and $g$ are treated
in fact separately. To prove a result with $g$ being either
H\"older or  Lipschitz continuous requires the rather
sophisticated regularization argument of case (3) and is therefore
of a different level of difficulty.

\noindent 2. The case of ``regular'' oblique derivative boundary
condition does not differ so much technically from the first case. The
assumption which says that $A(x)\gamma (x) = n(x)$ where, for any $x\in
\partial O$, $A(x)$ is a nonnegative symmetric matrix with a $C^2$
dependence in $x$, implies that $\gamma$ is a $C^1$ function
of $x$.

\noindent 3.  The case when $\gamma$ is only Lipschitz
continuous and when $g$ is only H\"{o}lder or Lipschitz
continuous, is technically very different as it is for the
comparison results (cf.  \cite{BaNeu}).  Here the only (known) way to
treat this case is through a non-trivial regularization argument
which we are going to use also.

Of course, in the proof below, we emphasize more the (difficult) third
case~:
the proofs are far easier in the two first ones.

Assumption {\bf (H2)} is a classical hypothesis in such type of
regularity result~: it is the same as the one which appears for
the interior regularity (cf. Ishii \& Lions\cite{il},
Barles\cite{Ba1}). It is worth pointing out, anyway, that the treatment
of the oblique derivative boundary condition does not lead to a
stronger assumption.

Concerning {\bf  (H3a)}, we recall that, for the interior regularity,
just the strong ellipticity ``in the gradient direction'' (cf.
(\ref{egd})) is needed as in the case of homogeneous boundary
condition. Unfortunately, in the case of
oblique derivatives boundary conditions, this natural assumption does
not seem to be enough or, at least, it has to be reformulated in a far
less natural way. Of course, all these conditions hold when a classical
uniform ellipticity property holds, like {\bf (H3b)} below.

One of the main examples we have in mind is the case of standard
quasilinear equations
      \begin{equation}\label{qle}
-{\rm Tr} [b(x,Du)D^2u]+H(x,u,Du)=0~~\mbox{in $O,$}
\end{equation}
where $b$ is a $n\times n$ matrix and $H$ a continuous function.
In this case, the assumptions {\bf  (H2)} and {\bf (H3a)} are
easily checkable.

{\bf  (H3a)} is equivalent to~: there exists $\lambda >0$ such
that, for any $x\in \overline O$, $p \in \R$
$$ b(x,p) \geq \lambda  \widehat{A^{-1}(x)p} \otimes
\widehat{A^{-1}(x)p} - o(1)Id \; ,$$
where, as in {\bf  (H3a)}, $o(1)$ is a function of $|p|$ which
converges to $0$ as $|p|\to +\infty$.

This assumption may seem restrictive, in particular the fact that the
constant $\lambda$ does not depend on $x$ and $p$; but, for general
equations, if $b$ satisfies the above statement with a strictly
positive $\lambda$ depending on $x$ and $p$, one can divide the
equation (i.e $b$ and $H$) by $\lambda (x,p)$ and the above property
becomes true with $\lambda = 1$.

We conclude these remarks about {\bf  (H3a)} by emphasizing the role of
the ``$o(1)$'' term and, for the sake of simplicity, we assume that
$A\equiv Id$. Without this term, {\bf  (H3a)} would be essentially
reduced to
$$ b(x,p) \geq \lambda  \hat{p} \otimes  \hat{p} \; ,$$
for any $x \in \overline O$ and $p \in \R-\{0\}$,  while, with this
term,  {\bf  (H3a)} is satisfied if
$$ b(x,p) \geq \lambda  \widehat{q (x,p)} \otimes  \widehat{q (x,p)}
\; ,$$
where $q$ is a continuous function such that $|p|^{-1}(q(x,p) - p) \to
0$ as $p \to \infty$, uniformly with respect to $x\in \overline O$.
This condition is not only a more general assumption on $b$ but it is
also far easier to check it. Of course, a similar remark can be made
for general $A$'s.

Now we turn to {\bf (H2)}. It is satisfied when

(i) $b$ is a bounded, continuous function of $x$ and $p$ and there
exists a modulus of continuity $\omega_{1}: \r^+ \to \r^+$ and a
function $\varpi : \r^+ \to \r^+$ such
that $\varpi(t) \to 0$ as $t \to + \infty$ and
$$ |b(x,p)-b(y,q)| \leq  \omega_{1}(|x-y|(1+ |p|+|q|)
+\varpi(|p| \wedge |q|)|p-q|)\; .$$ Moreover the uniform bound on
$b$ provides a $\omega_{2}$ with a linear growth.\par

(ii) The function $H$ satisfies~: for any $R>0$, there
exist positive constants $C_1^R$, $C_2^R$, $C_3^R$ such that, for any
$x,y \in
\overline O$, $-R\leq u,v \leq R$ and $p,q \in \R$,
$$ H(x,u,p)-H(y,v,q) \leq
C_1^R + C_2^R(|p|^2+|q|^2)+C_3^R|x-y|(|p|^3+|q|^3)\;  .$$
As we already mentioned it above, this assumption is classical (See
Ishii \& Lions\cite{il},
Barles\cite{Ba1}).

\noindent{\bf Proof of Theorem~\ref{theorregul}~: }
We are going to do the proof in two steps~: the first one consists in
proving the result
  when $a\equiv 0$ and contains the main arguments. In the second one,
we show how to handle the ``$a(x) u$'' term.
   According to the proof of the comparison result, it is clear that to
take in account such a $u$-term in regularity result is not immediate
and we are going to do it in a very indirect way, by adding an extra
variable.

\medskip
\noindent{\bf Step 1 : The $a\equiv 0$ case.}\newline
Since we are going to argue locally, we start with some notations.  For
every $x_0\in\overline
O$ and $r>0$, $B_{\overline O}(x_0,r):=B(x_0,r)\cap\overline O$
denotes the open ball in the topology of $\overline O$ while
$\overline{B}_{\overline
O}(x_0,r):=\overline{B}(x_0,r)\cap\overline O$ denotes the closed
ball in the topology of $\overline O$ and $\partial{B}_{\overline
O}(x_0,r):=\partial {B}(x_0,r)\cap\overline O.$ In a similar way
we define $B_{O}(x_0,r),$ $\overline{B}_{O}(x_0,r),$ and
$\partial{B}_{O}(x_0,r).$

We assume that $g\in C^{0,\beta}_{loc}(\partial O)$ for some
$0<\beta<1$, the case $g\in C^{0,1}_{loc}(\partial O)$ being treated in
a similar
way.

We are going to prove that, if we choose $\alpha = \beta$, for all
$x_0\in\overline O$, if $r$ is small enough, then there exists a
constant $C$ depending on the different data of the problem, such that,
for any $x \in B_{\overline O}(x_0,r)$, we have
\begin{equation}\label{holderest}
u(x_0)-u(x) \leq C|x-x_0|^\alpha\; .\end{equation} All together,
these inequalities give the answer provided that we control the
dependence of $r$ and $C$ in $x_0$, which will be the case.

The proof of this estimate is done in two steps whose arguments are the
same: the first step consists in proving that the result holds for
$\alpha$ small enough (depending on the local $L^\infty$ norm of $u$
and the data of the problem) and then that this property implies that
the result holds for $\alpha = \beta$.

We provide the main arguments of the proof of (\ref{holderest}) in the
case when the boundary condition plays a role, i.e. when $B(x_0, r)
\cap \partial O\neq \emptyset $; the other case is simpler and can be
treated by the methods of \cite{Ba1}.
In the common proof of the two steps, we therefore argue with some
$\alpha \leq \beta$. The following lemma is the key stone of the proof.

\begin{lemma}\label{keylemma} Assume that $B(x_0, r) \cap \partial
O\neq \emptyset $ and that $u$ is a bounded, continuous solution of
(\ref{np}) in $B_{\overline O}(x_0, 3r)$ with the oblique derivative
boundary condition on $B(x_0,3r)\cap\partial O$. Under the assumptions
of Theorem~\ref{theorregul} on $F$,
$\gamma$ and $g$, then there exists a constant  $C>0$, depending on
$F,\gamma,g$ and $||u||_{L^\infty(\overline B_{\overline O}(x_0, 3r))}$
such that for all $x\in B_{\overline O}(x_0,r) $ the
estimate \rec{holderest} holds.
\end{lemma}

{\noindent\bf Proof.}  In order to prove \rec{holderest}, we consider
the auxiliary function
$$ \Phi_0(x,y)=u(x)-u(y)-
\Theta_{0}(x,y)\; , $$ where the function $\Theta_{0}$ has the
following form $$  \Theta_{0}(x,y)= Ce^{-\tilde K(d(x)+d(y))}
[\psi_0(x,y)]^{\frac{\alpha}{2}}+ Le^{-\tilde K(d(x)+d(x_0)
}\psi_0(x,x_0) +\chi_0 (x,y)\; ,$$
where $\alpha \in (0,\beta]$ is a fixed constant, $C$, $L$, $\tilde K$
are some large constants to be chosen later on and where the continuous
functions $\psi_0(x,y)$, $\chi_0(x,y)$ satisfy the properties listed in
the lemma below. In order to point out the main
dependences in these functions but also to simplify the rather
technical estimates we have to make in the proofs, we introduce
the following notations which are used in all the sequel $$ X :=
\frac{x+y}{2}\; ,\; Y := x-y\; ,\; Z := d(x)-d(y)\; ,\; T := d(x)
+ d(y)\; .$$

\begin{lemma}\label{lemmaprel}
Under the assumptions of Theorem~\ref{theorregul}, for $\delta \geq 0$
small enough, there exist real-valued, continuous functions $\tilde
\psi_\delta(X,Y,T)$, $\tilde \chi_\delta
(X,Y,T,Z) $ defined respectively in $\R \times \R\times \r$ and
in $\R \times \R\times \r \times \r$, such that if we set $$
\psi_\delta(x,y):= \tilde
\psi_\delta(\frac{x+y}{2},x-y,d(x)+d(y))\;,$$ $$
\chi_\delta(x,y):= \tilde
\chi_\delta(\frac{x+y}{2},x-y,d(x)+d(y),d(x)-d(y))\; ,$$ the
following facts hold for some constant $K$ depending only on the
local $C^{0,1}$-norm of $\gamma$ and the $C^{0,\beta}$-norm of
$g$.

(i) For any $X,Y\in \R$, $T \in  \r$ and for $\delta \geq 0$ small
enough,
\begin{equation}\label{estpsio}
   K^{-1}|Y|^2\le \tilde \psi_\delta (X,Y,T) \le   K |Y|^2+K\delta\; ,
\end{equation}
\begin{equation}\label{estchio}
   -K|Z| \le \tilde \chi_\delta (X,Y,T,Z ) \le  K|Z|+K\delta^{\alpha}\; .
\end{equation}

(ii) When $\delta \to 0$, $\tilde \psi_\delta(X,Y,T) \to \tilde
\psi_0 (X,Y,T)$ and $\tilde \chi_\delta (X,Y,T,Z) \to \tilde
\chi_0 (X,Y,T,Z)$, uniformly on each compact subset of $\R \times
\R \times  \r$ and $\R \times
\R \times  \r  \times  \r$ respectively.

(iii) When $\delta>0$, the functions $\tilde \psi_\delta$, $\tilde
\chi_\delta$ are $C^2$ functions and
the following estimates hold for any $X,Y,Z,T$
\begin{eqnarray*}
&& \langle D^2_{{Y}{Y}}  \tilde\psi_\delta (X,Y,T) Y,Y\rangle
=2\tilde\psi_\delta +O(|Y|^3)+O(\delta) \;\mbox{as $Y\to
0\,,\,\delta\to 0$}\label{firsecderpsi}
\\
&& \langle D_{Y}\tilde \psi_\delta (X ,Y,T), Y\rangle
=2\tilde\psi_\delta +O(|Y|^3) +O(\delta) \;\mbox{as $Y\to
0\,,\,\delta\to 0$}\label{scalarprod}\\
&& |D_X \tilde \psi_\delta(X ,Y,T)|\le K|Y|^2\;,\; |D_{Y} \tilde
\psi_\delta (X,Y,T) | \leq K|Y|  \; \label{firstderpsi}\\ &&
|D_{{X}{X}}^2 \tilde \psi_\delta (X,Y,T)|\leq K|Y|
   \; , \;
|D_{{X}{Y}}^2 \tilde \psi_\delta (X ,Y,T)|\leq K |Y| \nonumber \\
&& |D_{YY}^2 \tilde \psi_\delta (X ,Y,T)|\leq K,\;,\; |D_T\tilde
\psi_\delta (X ,Y,T)|\leq K\delta\label{seconderpsi}\\
&& D_{TT}^2 \tilde \psi_\delta=D_{TX}^2 \tilde
\psi_\delta=D_{TY}^2 \tilde \psi_\delta=0\;,
\end{eqnarray*}
Moreover if  $\Lambda_\delta:= (|Y|^2+\delta^2)^{\frac{1}{2}}$
\begin{eqnarray*}
   &&  |D_{X} \tilde \chi_\delta (X ,Y,T,Z) | \leq
\Lambda^{\alpha-1}_\delta |Z|\; ,
   \\&&
|D_{Y}\tilde \chi_\delta (X ,Y,T,Z) | \leq
\Lambda^{\alpha-1}_\delta |Z|\; ,\;
\label{firstderchixy}
\\
&& |D_{Z}\tilde \chi_\delta (X ,Y,T,Z)|\le K\;,\;|D_{T}\tilde
\chi_\delta (X ,Y,T,Z)|\le K\delta ,\label{firstderchitz}\\&&
|D^2\tilde \chi_\delta (X ,Y,T,Z)|\leq
\Lambda_\delta^{\alpha-2} |Z| \; .\label{seconderchi}
\end{eqnarray*}
(iv) There exists a constant $\tilde K>0$ large enough (independent of
$C$
and $L$) such that, if we set
\begin{equation}\label{thetad1}
   \Theta_{\delta}(x,y)= Ce^{-\tilde K(d(x)+d(y))}
[\psi_\delta(x,y)]^{\frac{\alpha}{2}}+ Le^{-\tilde K(
d(x)+d(x_0)}\psi_\delta(x,x_0)+\chi_\delta(x,y)\; ,
\end{equation}
then, for $|x-y|$ small enough,  we have
\begin{eqnarray}
&& \langle D_x\Theta_{\delta}(x,y), \gamma (x) \rangle + g(x)
0~~\mbox{if
$x\in\partial O \,,$} \label{bctheta1l}\\
&&\langle -D_y\Theta_{\delta}(x,y), \gamma (y) \rangle+ g(y)
<0~~\mbox{if
$y\in\partial O \,.$}\label{bctheta2l}
\end{eqnarray}
\end{lemma}

\medskip

The proof of the key Lemma~\ref{lemmaprel} is postponed to
Subsection 3.1.  The first two properties in the point (iii) in
Lemma \ref{lemmaprel}   are going to play a central role in the
proof.

We continue with the {\bf proof of Lemma~\ref{keylemma}}.

We are going to show that, for a suitable choice of $L>0$, chosen large
enough
in order to localize, then for $C>0$ large enough, we have
\begin{equation}\label{maxneg}
M_{L,C}:=\max_{\overline{B}_{\overline  O}(x_0,r)\times
\overline{B}_{\overline  O}(x_0,r)} \Phi_0(x,y)\le 0
\end{equation}
Indeed if \rec{maxneg} holds then plugging $x=x_0$ and using the
estimates \rec{estpsio} and \rec{estchio} we get
$$ u(x_0)-u(y)\le K|x_0-y|^\alpha \; ,$$
for some constant $K>0$ depending on $\alpha$, $r,$
$||u||_{L^\infty(\overline{B}_{\overline  O}(x_0,r))}$, $F$ through
{\bf  (H2)-(H3a)}, the uniform  Lipschitz norm of $\gamma$ and  the
uniform H\"older norm of $g$ in $\overline{B}_{\overline
O}(x_0,r)$ and independent of $x_0$ (at least if $x_0$ remains in a
compact subset of $\overline O$), which is the desired regularity
result.

To prove (\ref{maxneg}), we first choose $L,C>0$ large enough in order
to have $\Phi_0(x,y) \leq 0$ for $x$ or $y$
on $\partial B (x_0,r) \cap \overline O$. This is possible since $u$ is
locally bounded on $\overline O$ and since
the conditions \rec{estpsio} and \rec{estchio} imply that
\begin{equation}\label{estpsio2}
   C_1|x-y|^2\le\psi_0(x,y)\le   C_2|x-y|^2 \; ,
\end{equation}
\begin{equation}\label{estchio2}
   -C_3|x-y| \le\chi_0(x,y)\le   C_4|x-y|\; .
\end{equation}
Of course, $L$ and $C$ depends on $r$.

 From now on, we fix such an $L$ and we argue by contradiction
assuming that, for  all $C>0$, $M_{L,C} >0$. Since $\Phi_0$ is a
continuous function, the maximum is
achieved at some $(\ox,\oy)\in \overline{B}_{\overline O}(x_0,r)
\times \overline{B}_{\overline O}(x_0,r)$ and we observe that, by
the choice of $L,C$, we may even assume that $\ox \in B_{\overline
O}(x_0,3r/4)$ and $\oy \in B_{\overline O}(x_0,r).$ Here we have
dropped the dependence of $\ox,\oy$ on $C$ for
simplicity of notations.

Two quantities are going to play a key role in the proof
$$ Q_1 := C |\ox-\oy|^\alpha\; ,$$
$$Q_2 :=  L|\ox-x_0|^2\; ,$$
(again we have dropped the dependence of $Q_1, Q_2$ in $C$ for the sake
of simplicity of notations).
The reason for that is the following: by using only the local
boundedness of $u$, we are only able to show that $Q_1, Q_2$ are
uniformly bounded when $C$ becomes very large while if we use the local
modulus of continuity of $u$, we can show that $Q_1, Q_2 \to 0$ as $C
\to +\infty$. The idea of the proof can therefore be described in the
following way: we first show that $u$ is locally in $C^{0, \alpha}$ for
$\alpha$ small enough with suitable estimates depending only on the
local $L^\infty$ norm of $u$ and on the data, and this is done by using
only the uniform boundedness od $Q_1, Q_2$. Then this first step
provides us with a local modulus of continuity for $u$ and we obtain
the full result using this time that $Q_1, Q_2 \to 0$ as $C \to +
\infty$.

As we just mention it, from the fact that $\Phi_0(\ox,\oy)>0$, using
classical arguments,
$Q_1, Q_2$ are bounded and, more precisely the following estimates
hold, in which $\bar K$ denotes the
constant $\left( e^{2\tilde K
r}(2||u||_\infty+||\chi_0||_\infty)\right)^{1/\alpha}$
\begin{eqnarray}\label{estxminusy}
Q_1 \leq \bar K ^\alpha & \hbox{or equivalently} &|\ox-\oy|\le \bar K
C^{-1/\alpha}\; , \\
&  L|\ox-x_0|^2\le \bar K & \,.\nonumber
\end{eqnarray}
  From the first estimates  \rec{estxminusy} it follows, in particular,
that
$|\ox-\oy| \to 0$ as $C\to + \infty.$ These estimates hold true for
any maximum point of the function $\Phi_0$ in
$\overline{B}_{\overline  O}(x_0,r) \times \overline{B}_{\overline
O}(x_0,r)$.

The function $\Theta_0$ is not (a priori) a smooth function and
therefore we cannot use directly viscosity solutions arguments;
this is why we have to consider the functions $\psi_\delta$ and
$\chi_\delta$ defined in Lemma~\ref{lemmaprel}.  Since
$\psi_\delta\to\psi_0$ and $\chi_\delta\to\chi_0$ as $\delta\to 0$
locally uniformly in $\overline O\times\overline O$, for all
$L,C>0$, there is $\delta_{C,L}>0$ such that for
$0<\delta\le\delta_{C,L}$, we have
\begin{equation}\label{maxneg2n}
\max_{\overline{B}_{\overline  O}(x_0,r)\times \overline{B}_{\overline
O}(x_0,r)}\,\left(u(x)-u(y)-\Theta_\delta(x,y)\right)> 0\,.
\end{equation}
Let $(x^\delta,y^\delta )$ be the maximum point of the function
$(x,y) \mapsto u(x)-u(y)-\Theta_\delta(x,y)$ in
$\overline{B}_{\overline O}(x_0,r)\times \overline{B}_{\overline
O}(x_0,r).$ Standard arguments show that, up to subsequence,
$(x^\delta,y^\delta)$ converges to a maximum point $(\ox,\oy)$ of
$\Phi_0$ as $\delta\to 0.$ Moreover we may suppose that
$x^\delta-y^\delta\ne 0.$ Indeed if for all $\delta>0$ we have
$x^\delta-y^\delta=0,$ then $\ox-\oy=0$ as well.  But in this case
we would have $\Phi_0(\bar x,\bar y) \le 0$ which is a
contradiction. Hence, we can assume without loss of generality
that $x^\delta-y^\delta$ remains bounded away from $0.$ For
simplicity of notations we now drop the dependence of
$(x^\delta,y^\delta)$ on $\delta$ as we already dropped it on $C.$
For $C$ large enough, we have $x,y \in B_{\overline O}(x_0,r)$.
Moreover, from Lemma~\ref{lemmaprel}, it follows that
\begin{eqnarray}  & & \langle D_x\Theta_{\delta}(x,y),\gamma (x)\rangle
+g(x) >0
~~\mbox{if $x\in\partial O$\, ,}\label{boundcondx}
\\
   && \langle - D_y \Theta_{\delta}(x,y), \gamma (x)\rangle +g(y)
<0~~\mbox{if
$y\in\partial O\,.$}\label{boundcondy}
\end{eqnarray}
Thus the viscosity inequalities associated to the equation $F=0$
hold for $u(x)$ and $u(y)$ whenever $x,y$ lie.

By the arguments of User's Guide \cite{cil}, for all $\eps>0$, there
exist $(p,B_1)\in \overline{J}^{2,+}u(x),$ $(q,B_2)\in
\overline {J}^{2,-}u(y)$ such that
$$ p=D_x\Theta_{\delta}(x,y) \quad , \quad
q=-D_y\Theta_{\delta}(x,y)\; ,$$
\begin{equation}\label{matrix}
-(\e^{-1}+||D^2\Theta_{\delta}(x,y)\vert\vert)Id\leq
\left(\begin{array}{cc} B_1 & 0 \\ 0 & -B_2 \end{array}\right)
   \leq  D^2  \Theta_{\delta}(x,y)+\e(D^2
   \Theta_{\delta}(x,y))^2 \, ,
\end{equation}
and
\begin{equation}\label{ineq}
F(x,u(x), p,B_1)\le 0 \quad , \quad F(y, u(y),q,B_2)\ge 0.
\end{equation}

We choose below $\e = \rho ||D^2\Theta_{\delta}(x,y)||^{-1}$ for $\rho$
small enough but fixed. Its size is determined in the proof below. Next
we need the following
lemma. We recall that $Y$ denotes $x-y$ and $\hat Y = \displaystyle
\frac{Y}{|Y|}$.
\begin{lemma}\label{estab}
If $\rho$ is small enough and if $B_1,B_2$ satisfy \rec{matrix}
then, for $|Y|$ small enough (i.e. for $C$ large enough), there is
  $K>0$ such that
\begin{eqnarray}\label{estimab}
\langle (B_1-B_2){\hat Y},{\hat Y}\rangle  & \leq & -
CK^{-1}\alpha(1-\alpha)
\frac{(\tilde\psi_\delta)^{\frac{\alpha}{2}}}{ |Y|^2}+
O(\delta)|Y|^{\alpha-4} \\ &&+ CK\vert Y\vert^{\alpha-1}+K
\,,\end{eqnarray} as $\delta \to 0.$
  Moreover $B_1-B_1\le \Ktu Id$ with $\Ktu $ given by
  $$
\Ktu =K(C\delta\vert Y\vert^{\alpha-2}+C\vert
Y\vert^{\alpha-1}+1)\,. $$
\end{lemma}

{\noindent\bf Proof of Lemma~\ref{estab}.}  By the regularity
properties of $\psi_\delta, \chi_\delta$ given in
Lemma~\ref{lemmaprel}, it is tedious (but easy) to check that all
the terms in $D^2\Theta_{\delta}(x,y)$ are estimated by $K+CK\vert
Y\vert^{\alpha-1} $ except perhaps the ones coming from the
derivation of the first term. More precisely, we have $$
D^2\Theta_\delta(x,y)=M_1+M_2+M_3+M_4 $$ where $$ M_1=Ce^{-\tilde
K(d(x)+d(y))}\left[\frac{\alpha}{2}(\psi_\delta)^{\alpha/2
-1}D^2\psi_\delta+
\frac{\alpha}{2}(\frac{\alpha}{2}-1)(\psi_\delta)^{\alpha/2
-2}D\psi_\delta\otimes D\psi_\delta\right]\,, $$ $$
M_2=C(\psi_\delta)^{\alpha/2}D^2(e^{-\tilde K(d(x)+d(y))})\,, $$
$$ M_3=C{\alpha}(\psi_\delta)^{\alpha/2-1}(D(e^{-\tilde
K(d(x)+d(y))})\otimes D\psi_\delta) \,, $$ $$ M_4=L D^2(e^{-\tilde
Kd(x)}\psi_\delta(x,x_0))+D^2\chi_\delta(x,y)\,. $$ But
we have also
\begin{eqnarray*}
D_x\psi_\delta (x,y)
&=&\frac{D_X\tilde\psi_\delta}{2}+D_Y\tilde\psi_\delta+
D_T\tilde\psi_\delta Dd(x) \,,
\\
D_y\psi_\delta (x,y)
&=&\frac{D_X\tilde\psi_\delta}{2}+D_Y\tilde\psi_\delta+
   D_T\tilde\psi_\delta Dd(y)  \,,\\ D^2_{xx}\psi_\delta (x,y)
&=&\frac{D^2_{XX}\tilde\psi_\delta}{4}+D^2_{YY}\tilde\psi_\delta+
D^2_{XY}\tilde\psi_\delta +   D_T\tilde\psi_\delta D^2d(x) \,,
\\
D^2_{yy}\psi_\delta (x,y)
&=&\frac{D^2_{XX}\tilde\psi_\delta}{4}+D^2_{YY}\tilde\psi_\delta-
D^2_{XY}\tilde\psi_\delta +  D_T\tilde\psi_\delta  D^2d(x) \,,
\\
D^2_{xy}\psi_\delta (x,y)
&=&\frac{D^2_{XX}\tilde\psi_\delta}{4}-D^2_{YY}\tilde\psi_\delta\,.
\end{eqnarray*}
and, taking in account the properties given in Lemma~\ref{lemmaprel},
it can be readily checked that $||M_2||\le CK|Y|^{\alpha},$ $||M_3||\le
   CK\alpha(|Y|^{\alpha-1}+\delta|Y|^{\alpha-2})$ and $||M_4||\le K \,.$

On the other hand, for all $\xi,\zeta \in\r^n$, we have
\begin{eqnarray}\label{estM1}
   \langle M_1 (\xi,\zeta),  (\xi,\zeta) \rangle & = & Ce^{-\tilde
K(d(x)+d(y))}\big[
\langle P_{1}
(\xi-\zeta),(\xi-\zeta)\rangle + \langle P_{2}
(\xi+\zeta),(\xi-\zeta)\rangle \nonumber \\
&& +
   \langle P_{3} (\xi+\zeta),(\xi+\zeta) + \langle P_4 (\xi,\zeta),
(\xi,\zeta) \rangle\big]
\end{eqnarray}
where $$P_{1}= D_{YY}^2 (\tilde{\psi}_\delta)^{\alpha/2} (X,Y,T)\;
,\;P_{2}= 2D_{XY}^2  (\tilde{\psi}_\delta)^{\alpha/2}  (X,Y,T)\;
,\;P_{3}= D_{XX}^2 (\tilde{\psi}_\delta)^{\alpha/2} (X,Y,T)\; ,$$
and $P_4$ is the matrix involving all the terms
$D^2_{XT}(\tilde{\psi}_\delta)^{\alpha/2} ,
D^2_{YT}(\tilde{\psi}_\delta)^{\alpha/2},
D^2_{TT}(\tilde{\psi}_\delta)^{\alpha/2}$ together with the $Dd$
and $D^2d$ derivatives. One can easily check that $P_4$ is
estimated by $O(\delta)(\tilde{\psi}_\delta)^{(\alpha -
4)/2}$. On the other hand, Lemma~\ref{lemmaprel} implies $$ |P_1|
= O( (\tilde{\psi}_\delta)^{(\alpha-2)/2})\; ,$$ $$ |P_2|  = O(
(\tilde{\psi}_\delta)^{(\alpha-1)/2})\; ,$$ $$ |P_3|= O(
(\tilde{\psi}_\delta)^{(\alpha-1)/2})\; .$$

Choosing $\xi=\zeta$ in \rec{estM1}, we first deduce that
$$ M_1\le C\left( O((\tilde{\psi}_\delta)^{(\alpha-1)/2})
+O(\delta)(\tilde{\psi}_\delta)^{(\alpha - 4)/2})\right) Id\;. $$
We next choose   $\xi=-\zeta={\hat Y}.$ According to the two first
properties in the point (iii) of Lemma~\ref{lemmaprel} and taking
in account the above estimate on $P_4$, we get
\begin{eqnarray*}
\langle M_1 (\hat Y,-\hat Y), (\hat Y,-\hat Y) \rangle & \le & 4
Ce^{-2\tilde K(d(x)+d(y))}
\frac{\alpha}{2}\frac{(\tilde\psi_\delta)^{\alpha/2}}{|Y|^2}
\left[
   \langle D_{YY}^2\tilde{\psi}_\delta Y,Y\rangle+
   (\frac{\alpha}{2}-1)\frac{ \langle D_Y\tilde\psi_\delta,
Y \rangle^2}{\tilde\psi_\delta}\right]\\
&& + C\langle P_4 (\hat Y,-\hat Y), (\hat Y,-\hat Y) \rangle
\\&\leq &
CK^{-1} \alpha(\alpha-1)\tilde\psi_\delta^{(\alpha-2)/
2}+K
(\tilde\psi_\delta)^{\alpha/2}+O(\delta)\tilde\psi_\delta^{(\alpha-4)/
2}\,.
   \end{eqnarray*}

By combining the above estimates, we obtain
\begin{equation}\label{estb1b2l}
\langle B_1-B_2 \hat Y,\hat Y\rangle \le CK^{-1}\alpha(\alpha-1)
|Y|^{\alpha-2}+ C O(\delta)|Y|^{\alpha-4}+ CK\vert
Y\vert^{\alpha-1}+K \,.
\end{equation}
And the final upper estimate on $B_1-B_2$ follows from the
estimates on $M_i$ for $i=1,2,3,4$.~~\hfill$\square$

\medskip We continue with the {\bf proof of Lemma~\ref{keylemma}}~: we
estimate $|p-q|,$ $|p|,|q|$ and $||B_1||, ||B_2||$.
For some $K>0$, we have
$$ |p|,|q|\ge  CK^{-1}\alpha
|Y|^{\alpha-1} - O(\delta)|Y|^{\alpha-3}- K\; ,$$
$$ |p|,|q| \le
CK\alpha|Y|^{\alpha-1}+ O(\delta)|Y|^{\alpha-1}+CK|Y|^\alpha+K \; ,$$
$$
   |p-q| \le KC\alpha|Y|^{\alpha}+O(|x-x_0|)+o_\delta(1),~~\mbox{as
$|Y|\to
   0,~\delta\to 0\,,$}
   $$
$$ ||B_1||, ||B_2||\le
K(1+\displaystyle\frac{1}{O(\rho)})\left[C\alpha |Y|^{{\alpha}
-2}+O(\delta)|Y|^{\alpha-4}
   + 1\right]\,. $$
At this point, it is worth noticing that we are going to let first
$\delta$ tends to $0$ for fixed $C$ and we recall that, since we assume
that $M_{L,C}>0$, $Y$ does not converge to $0$  when $\delta$ tends to
$0$ for fixed $C$. The first consequence of this fact is that the term
$O(\delta)|Y|^{\alpha-3}$ is playing no role in the lower estimate of
$|p|, |q|$ since we can choose $\delta$ as small as necessary and, by
the above estimates, we have $|p|, |q| \to + \infty$ as $C \to +\infty$

We are going to use {\bf  (H2)-(H3a)} with
$R=||u||_{L^\infty(B_{\overline O}(x_{0},r))}\,.$ We drop the
dependence in $R$ in the coefficients and modulus which appear in these
assumptions.  We subtract the two inequalities \rec{ineq}
and write the difference in the following way
   \begin{multline}\label{diff}
F(x,u(x),p,B_1)- F(x,u(x),p, B_2+\Ktu Id)\\
\leq F(y,u(y),q,B_2) - F(x,u(x),p, B_2+\Ktu
Id)\; ,
\end{multline}
and, using the fact that $B_1-B_2 \leq \Ktu  Id$, we apply
{\bf  (H3a)} to the left-hand side and {\bf (H2)} to the
right-hand side of \rec{diff}. Recalling also that $|p|, |q| \to +
\infty$ as $C \to +\infty$, this yields
\begin{multline}\label{eqnfin}
\lambda\mbox{Tr}[(B_2-B_1+\Ktu Id)(\widehat{A^{-1}(x)p}\otimes
\widehat{A^{-1}(x)p})]  + o(1) ||B_2-B_1+\Ktu Id ||  \\
\le \omega_1 (|x-y|(1+|p|+|q|) +\varpi(|p|\wedge|q|)|p-q|) ||B_2||
+ \omega_{2}(\Ktu)
  + C_1 \\
+C_2 (|p|^2+|q|^2)+C_3|x-y|(|p|^3+|q|^3)
   \end{multline}

Now we use the following result which is a consequence of the
construction of the test-function and whose proof is given at the end
of Subsection 3.1.
\begin{lemma}\label{lemmaAY}
   We have
   $$
\widehat{A^{-1}(x)p}={\hat Y}+o_Y(1)+o_\delta(1)~~~\mbox{as
$|Y|\to 0\,,\delta\to 0\,.$}$$
\end{lemma}

We want to point out that, in the above lemma, $|Y|\to 0$ is in fact
equivalent to $C$ going to infinity.
\par\medskip
We come back to (\ref{eqnfin}): by  Lemma \ref{lemmaAY} and recalling
also that $|p|, |q| \to + \infty$ as $C \to +\infty$, we get
\begin{multline*}
\mbox{Tr}[(B_2-B_1+\Ktu Id)(\widehat{A^{-1}(x)p}\otimes
\widehat{A^{-1}(x)p})] \\ \ge
\langle (B_2-B_1){\hat Y},{\hat Y}\rangle+\Ktu
-||B_2-B_1||(o_Y(1)+o_\delta(1)) \;.
\end{multline*}
Moreover by using the estimates on $B_1,B_2$ we have
  {\bf  \begin{equation}\label{estb1b2n}
||B_2-B_1||\le
K\left[1+\displaystyle\frac{1}{O(\rho)}\right]\left[C\alpha |Y|^{{\alpha
}-2}+O(\delta)|Y|^{\alpha-4}
   +1\right]\;.
\end{equation}}
As we already pointed out above, we are going to let first $\delta$
tends to $0$ for fixed $C$ and, since we assume that $M_{L,C}>0$, $Y$
does not converge to $0$  when $\delta$ tends to $0$ for fixed $C$.
Hence in the estimates below, we are going to replace the terms which
converge to $0$ as $\delta \to 0$ by $o_\delta(1)$. On the other hand,
as we already mention it above, $C$ going to infinity is equivalent to
$Y$ going to $0$ and when $C$ is going to infinity, $p,q$ are also
going to infinity; we can therefore incorporate the $o(1)$-term coming
from {\bf (H3a)} in the $o_Y(1)$ term.

Therefore, by combining \rec{estb1b2n}, the estimates on
$||B_2||$, $|p|,|q|$ and $|p-q|$, and Lemma \ref{lemmaAY} we are
lead to
\begin{multline*}
\mbox{Tr}[(B_2-B_1+\Ktu Id)(\widehat{A^{-1}(x)p}\otimes
\widehat{A^{-1}(x)p})] \\
\ge    C K^{-1} \alpha(1-\alpha)|Y|^{\alpha-2} + o_\delta(1)-
K \\
-(CK\alpha|Y|^{\alpha-2}+o_\delta(1)+K)(o_Y(1)+o_\delta(1))  \; .
\end{multline*}

\medskip

On the other hand, for the right-hand side of \rec{eqnfin}, we first
look at the $\omega_1$ term. By tedious but straightforward
computations, we have
$$ |x-y|(1+|p|+|q|) +\varpi(|p|\wedge|q|)|p-q| = K\alpha Q_1 + K
\varpi(|p|\wedge|q|) Q_2^{1/2} + o_Y(1) + o_\delta(1)\; ,$$
since $O(|x-x_0|)$ is like $Q_2^{1/2}$. This estimate is emphasizing
the role of $Q_1, Q_2$ and the necessity of having the $\varpi$ term.

The complete estimate of the right hand side of  \rec{eqnfin} is
\begin{multline*}
K \omega_1 \left( K\alpha Q_1 + K \varpi(|p|\wedge|q|) Q_2^{1/2} +
o_Y(1) + o_\delta(1) \right) C\alpha|Y|^{\alpha - 2}  \\ + K
C^2\alpha^2|Y|^{2\alpha-2} + C^3 \alpha^3|Y|^{3\alpha-2} + K+
o_Y(1)+o_\delta(1) \; ,
\end{multline*}
where we (partially) use the fact that $Q_1 = C|Y|^\alpha$  is bounded
for $C$ large enough.

  By dividing all the above inequalities  by the (very large) term $C
\alpha |Y|^{\alpha-2}$, we obtain the following
(almost) final estimate
\begin{multline*}
\lambda K^{-1} \leq
K \omega_1\left( K\alpha Q_1 + K \varpi(|p|\wedge|q|) Q_2^{1/2} +
o_Y(1) + o_\delta(1)\right)
  + K \alpha Q_1 + K \alpha^2 Q_1^2 \\ +  o_Y(1)+o_\delta(1)\; .$$
\end{multline*}
And by using the fact that $|p|, |q| \to + \infty $ as $C$ tends to $+
\infty$, this yields
\begin{eqnarray}\label{finalestim}
\lambda K^{-1}& \leq& K \omega_1\left( K\alpha Q_1 + o_Y(1)
Q_2^{1/2} + o_Y(1) + o_\delta(1)\right)\nonumber\\ & &
  + K \alpha Q_1 + K \alpha^2 Q_1^2 +  o_Y(1)+o_\delta(1)\;
.\end{eqnarray}

Using this last estimate, the conclusions of the two steps we
mention at the beginning of the proof follow rather easily.

On one hand, by using the uniform control on $Q_1, Q_2$, we can
choose $\alpha$ small enough (depending only on the local
$L^\infty$ norm of $u$ and the data) in order to have $$ \lambda
K^{-1} \geq  \frac{3}{2}\left(K \omega_1( K\alpha Q_1) + K \alpha
Q_1 + K \alpha^2 Q_1^2\right) > K\omega_1\left( K\alpha
Q_1)\right)
  + K \alpha Q_1 + K \alpha^2 Q_1^2 \; .$$
With this choice, it is clear that the above inequality cannot holds
for $\delta$ small and  $C$ large enough (depending again only on the
local $L^\infty$ norm of $u$ and the data) and the local $C^{0,\alpha}$
estimate is proved for small enough $ \alpha$.

On the other hand, repeating this proof for any $x_0 \in \overline
B_{\overline O}(x_0, 2r)$, this $C^{0,\alpha}$ property provides us
with a modulus of continuity in $B_{\overline O}(x_0, r)$ (which
depends only on the $L^\infty$ norm of $u$ in $B_{\overline O}(x_0,
3r)$ and the data), and in the above estimate, for any $\alpha \leq
\beta$, we can use the fact that $Q_1, Q_2 \to 0$ as $C\to + \infty$.
Arguing as above, we obtain the $C^{0,\alpha}$ estimate for any $
\alpha \leq \beta$. And the proof of the Step 1 is complete.
\hfill\sn \par

\medskip
\noindent{\bf Step 2 : how to handle the ``$a(x) u$''
term}\newline We are going to introduce an extra variable to
reduce this case to the previous one. More precisely we consider
the function $v \colon \overline O \times \r \to \r$ defined by
$$v(x,y) = \exp(ky)u(x)\, ,$$
where $k>0$ is a large constant to be chosen.

This new function is a solution of
$$ -\frac{\partial^2 v}{\partial
y^2}+ \exp(ky)F(x,\exp(- ky) v ,\exp(- ky)D_x v ,\exp(- ky) D_{xx}^2v)
+ k^2 v =0\;\mbox{in $O\times \r\, ,$} $$
and
$$
\langle \gamma (x), D_x v \rangle+ \frac{1}{k} a(x) D_y v +
\exp(ky)g(x) =0\quad\mbox{on
$\partial O\times \r\,.$} $$

We first remark that since we are going to argue in a neighborhood of
the point $(x_0,0)$, the exponential terms are not a problem to check
the assumptions on either the equation or the boundary condition. The
only difficulty comes from {\bf (H3a)} since we have changed the
boundary, $\gamma$ and therefore $A$ is not given anymore.

In order to define the new matrix $A$, denoted below by $\tilde A_k$
since it depends on $k$, we first set $\tilde n := (n,0)$, the exterior
unit normal vector to $\partial O\times \r $ and $\tilde \gamma
:=(\gamma, k^{-1}a)$. In fact, because of the form of {\bf (H3a)}, it
is more convenient to define $\tilde A_k^{-1}$ and we do it by setting
$$ \tilde A_k^{-1}(x) :=
\left(   \begin{array}{cc} A^{-1}(x) &  k^{-1}a(x)n(x)\\
k^{-1}a(x)n^T(x) & 1\end{array} \right)\; ,$$ where $n^T$ is the
transpose of the column vector $n$.

An easy computation shows that $ \tilde A_k^{-1}(x) \tilde n (x)=\tilde
\gamma (x)$.
  Moreover, if $P=(p,p_{n+1})\in\r^{n+1}$, we have
$$ \langle \tilde A_k^{-1}(x) P,P\rangle =  \langle A^{-1} (x) p, p
\rangle + 2 k^{-1} a(x) \langle n (x) , p \rangle p_{n+1} + |p_{n+1}
|^2\; ,$$
and applying Cauchy-Schwarz inequality to the second term of the
right-hand side, it is straightforward to show that, for $k$ large
enough,  $ \tilde A_k^{-1}(x)$ is still a definite positive matrix.
Finally
we consider the operator $\tilde F\colon\overline
O\times\r\times\r\times\r^{n+1} \times {\cal{S}}^{n+1}\to\r$
defined by
$$
\tilde F(x,y,v,P,M)=-M_{n+1,n+1}
+ \exp(ky) F(x,\exp(- ky) v ,\exp(- ky)p ,\exp(- ky) \tilde M) + k^2 v
\,,
$$
where $\tilde M$ is the $n\times n$ symmetric matrix obtained from $M$
by removing the last
column and row.

Next we claim that the operator $\tilde F$ satisfies {\bf (H3a)} with $
\tilde A_k^{-1}$ but by replacing the $o_P(1)$ term which tends to $0$
as $|P| \to \infty$ by a $o_P(1) + o_k(1)$ where $o_k(1)\to 0$ as $k\to
+\infty\,.$

Indeed, we first observe that, since the new solution $v$ is Lipschitz
continuous w.r.t $y$ and since $y$ is a variable which corresponds to a
tangent direction to the boundary, the property ``$|P| \to + \infty$''
is equivalent in fact to ``$|p| \to + \infty$'' because $p_{n+1}$
remains bounded.

On an other hand, by using {\bf (H3a)} for $F$, one can easily see that
the checking of our property reduces to show that,
for all $R>0$, $x \in \overline O$, $P=(p,p_{n+1})$ and $
\zeta=(\xi,\eta)\in \r^{n+1}$, we have
$$ \lambda_R \langle A^{-1}(x) p, \xi\rangle^2 + |\eta|^2 \geq
\tilde\lambda_R  \langle \tilde A_k^{-1}(x) P, \zeta\rangle^2 -
o_k(1)\; ,$$
for some constant $\tilde\lambda_R >0$. Because of the particular form
of $\tilde A_k^{-1}(x)$, this property is obvious for ``$k=+\infty$"
and, of course, this implies that it is also satisfied for $k$ large
within a $o_k(1)$-term.

We finally observe that the $C^{0,1}_{loc}$ norm of $\tilde A_k^{-1}$
does not depend on $k$ if we choose it large enough.
In order to conclude, we just remark that the proof of the first step
still works if the term $o_k(1)$ is small enough and the proof is
complete.\hfill\sn

\medskip

\begin{remark} {\rm We remark that, under further regularity
assumptions on $O$ and the coefficients appearing in the boundary
condition \rec{oblder}, it is possible to handle the $a(x)u$ term
without adding an extra variable but by using another change
of variable. More precisely, let us suppose that the following
assumption
holds

\smallskip

\noindent {\bf (H5)} There exists a $C^2$-function $\chi :
\overline O \to \r$ such that $\displaystyle \frac{\partial
\chi}{\partial\gamma}= a(x)\quad \hbox{on  }\partial O\; .$

\smallskip

\noindent Then the  function $v$ defined by $v(x)=e^{\chi (x) }  u(x)$,
is a solution of a modified equation in $O$ (but still satisfying
{\bf (H2)-(H3a)}) with the boundary condition
$$ \frac{\partial v}{\partial\gamma}+  e^{\chi (x) } g(x) =0\quad
\hbox{on  }\partial O\; ,$$
  and we can apply the proof of the Step 1 of Theorem~\ref{theorregul}
to $v$.

Assumption {\bf (H5)} holds for example in the following case~: if
$O$ is a $C^{2,\beta}$ domain and $\gamma, a$ are $C^{1,\beta}$
function, then the existence of $\chi$ is given by Theorem~7.4 (p.
539) in Lieberman \& Trudinger\cite{LT}. Indeed, to build $\chi$,
one can solve the Laplace Equation in $O$ together with the
oblique derivative boundary condition.

Moreover, if $O$ is a $C^{3}$ domain and $\gamma, a$ are $C^{2}$, then
one
can just take $$\displaystyle \chi(x)=\frac{a(x) d(x)}{\langle
n(x),\gamma(x)\rangle}$$ in
$O$ where $a$, $\gamma$ and $n$ denotes here suitable extensions
of these functions to $\overline O$.}
\end{remark}
\medskip

\subsection{The case of nonlinear boundary conditions}

In this subsection, we consider the case of nonlinear boundary
conditions of the form
\begin{equation}\label{nlbc}
G(x,u,Du)=0\quad\mbox{in  $\partial O\;,$}
\end{equation}
where $G\colon\partial O\times\r\times\r^n\to \r$ is a continuous
function, satisfying the conditions {\bf(G1)} and {\bf(G2)}.

In this section, we use the following assumptions on $F$ and $G$.
\begin{enumerate}
\item[{\bf (H3b)}]~~{\bf (Uniform ellipticity)} For any $R>0$, there is
$\lambda_R > 0$ such that, for all $x\in\overline O$, $-R\leq u \leq
R$, $p \in \R$ and $M,N\in {\cal{S}}^n$ such that $M\le N$, we have $$
F(x,u,p,M)-F(x,u,p,N) \geq \lambda_R\mbox{Tr} (N-M)\;.$$
\end{enumerate}
\begin{enumerate}
\item[{\bf (G3)}]  For all $R>0$ and $M>0$ there is    $K_{R,M}>0$ such
that
\begin{equation}
|\langle\frac{\partial G}{\partial p}(x,u,p),p\rangle-G(x,u,p)|\le
K_{R,M}\,,
\end{equation}
\hspace{-\leftmargin}for all $x\in\partial O$ and for all $p\in\r^n$,
$|p|\ge M,$
$|u|\le R$\,.
    \item[{\bf (G4)}] There is a function $G_\infty\colon\partial
O\times\r\times\r^n\to\r$ such that
\begin{equation}
\frac{1}{\lambda}G(x,u,\lambda p)\to
G_\infty(x,u,p)\,\quad\mbox{as $\lambda\to\infty\,.$}
\end{equation}
\hspace{-\leftmargin}locally uniformly in $(x,u,p).$
\end{enumerate}

Before providing our result, we want to point out that the $G_\infty$
appearing in {\bf (G4)} is homogeneous of degree 1 and satisfies {\bf
(G1)} and {\bf (G2)}.

Our result is the
\begin{theor}\label{theorregulbis}
Assume {\bf (H1)-(H2)-(H3b)} and {\bf(G1)-(G4)}.  Then every
bounded continuous solution $u$ of \rec{np} is in
$C^{0,\alpha}_{loc}(\overline O)$ for any $0<\alpha<1$. Moreover
the $C^{0,\alpha}_{loc}$--norms of $u$ depend only on $O$, $F$,
$G$,   through the constants and functions appearing in {\bf
(H2)-(H3b)}, and in {\bf(G1)-(G4)}, the local $C^2$--norm of the
distance function of the boundary including the modulus of
continuity of $D^2 d$.
\end{theor}

{\noindent\bf Proof of Theorem~\ref{theorregulbis}.} We are going to do
the proof in three steps~: in the first one we prove the result in the
case when $G$ is independent of $u$ and homogeneous of degree 1 in $p$,
then in the second step, we remove the homogeneity restriction and
finally, in step 3, we use the method of the second step of the proof
of Theorem~\ref{theorregul} to deal with the dependence in $u$.

\smallskip
\noindent{\bf Step 1~: The case when $G$ is independent of $u$ and
homogeneous of degree 1 in $p$}\newline
The proof is similar to the one of Theorem~\ref{theorregul} and we just
outline the main differences. Again we treat only the case when the
boundary plays a role.

Since the boundary $\partial O$ is $C^2$, by making a suitable change
of variables, we can assume without loss of generality that the
boundary is flat and more precisely that $O\cap B(x_0, 3r)\subset \{x_n
0\}$ and $\partial O \cap B(x_0, 3r) \subset \{x_n = 0\}$.  It is
worth noticing that the assumptions made on $F$ and $G$ are preserved
by such a change.  In order to keep simple notations, we still denote
by $F$, $G$ the functions arising in the equation and in the boundary
condition in the domain with flat boundary.

We have to prove the following lemma which is the key stone of the
proof.
\begin{lemma}\label{flatbound} Assume that $B(x_0,r)\cap\{x_n >0\} \neq
\emptyset $ and that $u$ is a
continuous solution of (\ref{np}) in $B(x_0, 3r)\cap\{x_n > 0\} $
with nonlinear  boundary condition $G=0$ on $B(x_0, 3r)\cap\{x_n =
0\} $. Under the assumptions of Theorem~\ref{theorregulbis} on
$F$, $G$,  there exists a constant  $C>0$ depending on
$F,G,||u||_{L^\infty(B_{\overline O}(x_0, 3r))}$ such that, for any
$x\in B(x_0,r)\cap\{x_n \geq 0\}$,
the estimate \rec{holderest} holds.
\end{lemma}

{\noindent\bf Proof of Lemma \ref{flatbound} .}   In order to
prove \rec{holderest}, we  consider the auxiliary function $$
\Phi_0(x,y)=u(x)-u(y)- \Theta_{0}(x,y)\; , $$ where the function
$\Theta_{0}$ has in this case the following form $$
\Theta_{0}(x,y)= Ce^{-\tilde K(x_n+y_n)}
[\psi_0(x,y)]^{\frac{\alpha}{2}}+ L \phi_0(x,x_0)  \; ,$$ for some
large constants $C$, $L$, $\tilde K$ to be chosen later on and
where the continuous functions $\psi_0(x,y)$, $\phi_0(x,y)$
satisfy the properties listed in the following lemma in which we use
the notations
$$ X:= \frac{x+y}{2}\; ,\; Y:= x-y \; .$$

We want also to point out that the parameters $\delta$ and $\eta$ we
introduce in this lemma play completely different roles in the proof,
the role of $\delta$ being far more important than the role of $\eta$
which is a small but fixed parameter; this is why we choose to drop the
dependence in $\eta$ of the functions $\psi_\delta$, $\tilde
\psi_\delta$ below.

\begin{lemma}\label{lemmaprelbis}
Under the assumptions of Theorem~\ref{theorregulbis}, there is a
function $\phi_0\in C^2(\r^{2n})$ and,  for $\delta \geq 0$ and
$\eta>0$ small enough, there exists a real-valued, continuous function
$\tilde \psi_\delta(X,Y) $ defined in $\R \times \R$ such that if
we set $$ \psi_\delta(x,y):= \tilde \psi_\delta(\frac{x+y}{2},x-y)
$$ the following facts hold

(i) There exists a constant $K>0$ such that, for $\delta \geq 0$ small
enough,
\begin{equation}\label{estpsiobis}
   K^{-1}|Y|^2\le \tilde \psi_\delta (X,Y) \le   K |Y|^2+K\delta\; ,
\end{equation}
for any $Y=(Y_1, \ldots, Y_n) \in \R$.

(ii) When $\delta \to 0$, $\tilde \psi_\delta(X,Y) \to \tilde \psi_0
(X,Y)$  uniformly on each compact subset of $\R \times \R$.

(iii) When $\delta>0$, the function $\tilde \psi_\delta$ is $C^2$
   and the following estimates are valid for some constant $K>0$
$$ |D_X \tilde \psi_\delta (X,Y) | \leq K|Y|^2+K\delta\; ,
   \; |D_Y \tilde \psi_\delta (X,Y) | \leq K|Y|\; ,$$
   $$ \langle D_{Y}\tilde \psi_\delta (X,Y), Y\rangle =
2\tilde\psi_\delta(X,Y)+O(\eta) O(|Y|^2) +O(\delta)\; ,$$
$$  \langle D^2_{YY}\tilde \psi_\delta (X,Y) Y,Y\rangle =
2\tilde\psi_\delta(X,Y)+O(\eta) O(|Y|^2) +O(\delta)$$
as $Y\to 0,$ $\eta\to 0$ and  $\delta\to 0\,.$
   $$|D_{XX}^2 \tilde \psi_\delta (X,Y)|\leq \frac{K}{\eta}|Y|(\delta^2
+|Y|^2)^{1/2}\; , \;
|D_{XY}^2 \tilde \psi_\delta (X,Y)|\leq K |Y|  \; , \; |D_{YY}^2
\tilde \psi_\delta (X,Y)|\leq K.$$
  (iv)
\begin{equation}\label{estphio}
   K^{-1}|x-y|^4\le \phi_0 (x,y) \le   K
   |x-y|^4\, ,
   \end{equation}
(v) There exists $\tilde
K>0$ large enough (independent of $C$ and $L$) such that, if we
set
\begin{equation}\label{thetad2}
   \Theta_{\delta}(x,y)= Ce^{-\tilde K(x_n+y_n)}
[\psi_\delta(x,y)]^{\frac{\alpha}{2}}+ L \phi_0(x,x_0) \; ,
\end{equation}
then, for $|x-y|$ small enough,  we have
\begin{eqnarray}
&& G(x,D_x\Theta_{\delta}(x,y)) >0~~\mbox{if $x_n=0 \, ,$}
\label{bctheta1}\\
&& G(y, -D_y\Theta_{\delta}(x,y)) <0~~\mbox{if $y_n=0
\,.$}\label{bctheta2}
\end{eqnarray}
\end{lemma}

The proof of the key Lemma~\ref{lemmaprelbis} is postponed to the
Subsection 3.2. To prove (\ref{maxneg}), we first choose $L,C>0$
large enough in order to have $\Phi_0(x,y) \leq 0$ for $x$ or $y$
on $\partial B (x_0,r) \cap \overline O$. This is possible since
because of the conditions  \rec{estpsiobis} and
\rec{estphio}. Of course, $L, C$ depends on $r$.

As in the proof of Theorem \ref{theorregul}, we fix such an $L$
and we argue by contradiction assuming that, for all  $\alpha\in
(0,1)$ and $C>0$, $M_{L,C} > 0$.  Since $\Phi_0$ is a continuous
function, the maximum is achieved at some $(\ox,\oy)
\in\overline{B}_{\overline O}(x_0,r) \times
\overline{B}_{\overline O}(x_0,r)$ and we observe that, by the
choice of $L,C$, we may even assume that $\ox \in B_{\overline
O}(x_0,3r/4)$ and $\oy \in B_{\overline O}(x_0,r).$ Here we have
dropped the dependence of $\ox,\oy$ on $C$ for simplicity of
notations.

We use here
$ Q_1: = C|\ox-\oy|^\alpha$ and $Q_2 := L|\ox-x_0|^4$. From the fact
that $\Phi_0(\ox,\oy)>0$, using classical arguments,
the following estimates follow, in which $\bar K$ is the constant
   $\left( 2e^{2\tilde K r}||u||_\infty\right)^{1/\alpha}$
\begin{eqnarray}\label{estxminusybis}
Q_1 \leq \bar K^\alpha &\hbox{or equivalently} &|\ox-\oy|\le \bar K
C^{-1/\alpha} \\
&  Q_2 \le \bar K &  \,.\nonumber
\end{eqnarray}

    These estimates hold true for
    any maximum point of the function $\Phi_0$ in
     $\overline{B}_{\overline  O}(x_0,r) \times \overline{B}_{\overline
O}(x_0,r)$.

Since the function $\Phi_0$ is not (a priori) a smooth
function  we have to consider the functions $\psi_\delta$  defined in
Lemma~\ref{lemmaprelbis}.
    The following property holds~: for all $L,C>0$, there is
$\delta_{C,L}>0$ such that for $0<\delta\le\delta_{C,L}$ and $\tilde
K$ large, we have
\begin{equation}\label{maxneg2}
\max_{\overline{B}_{\overline  O}(x_0,r)\times \overline{B}_{\overline
O}(x_0,r)}\,\left(u(x)-u(y)-\Theta_\delta(x,y)\right)> 0\,.
\end{equation}

Let $(x^\delta,y^\delta )$ be the maximum point of
$u(x)-u(y)-\Theta_\delta(x,y)$ in $\overline{B}_{\overline
O}(x_0,r)\times \overline{B}_{\overline O}(x_0,r).$ Standard
arguments show that, up to subsequence, $(x^\delta,y^\delta)$
converges to a maximum point $(\ox,\oy)$ of $\Phi_0$ as $\delta\to
0.$ Again we may suppose that $x^\delta-y^\delta\ne 0.$   For
simplicity of notations we drop the dependence of
$(x^\delta,y^\delta)$ on $\delta$ as we already dropped it on $C.$

   For $C$ large enough, we have $x,y \in B_{\overline O}(x_0,r)$.
Moreover, from Lemma~\ref{lemmaprelbis}, it follows that
\begin{eqnarray}  & & G(x, D_x\Theta_{\delta}(x,y)) >0 ~~\mbox{if
$x_n=0 $,}\label{boundcondxbis}
\\
   && G(y,-D_y\Theta_{\delta}(x,y))  <0~~\mbox{if $y_n=0
\,.$}\label{boundcondybis}
\end{eqnarray}
Thus the viscosity inequalities associated to the equation $F=0$
hold for $u(x)$ and $u(y)$ whenever $x,y$ lie.

By the arguments of User's Guide \cite{cil}, for all $\eps>0$, there
exist $(p,B_1)\in \overline{J}^{2,+}u(x),$ $(q,B_2)\in
\overline {J}^{2,-}u(y)$ such that
$$ p=D_x\Theta_{\delta}(x,y) \quad , \quad
q=-D_y\Theta_{\delta}(x,y)\; ,$$
\begin{equation}\label{matrixbis}
-(\e^{-1}+||D^2\Theta_{\delta}(x,y)\vert\vert)Id\leq
\left(\begin{array}{cc} B_1 & 0 \\ 0 & -B_2 \end{array}\right)
   \leq  D^2  \Theta_{\delta}(x,y)+\e(D^2
   \Theta_{\delta}(x,y))^2 \, ,
\end{equation}
and
\begin{equation}\label{ineq2}
F(x,u(x), p,B_1)\le 0 \quad , \quad F(y, u(y),q,B_2)\ge 0.
\end{equation}

We choose below $\e = \rho ||D^2\Theta_{\delta}(x,y)||^{-1}$ for
$\rho$ small enough but fixed. Its size is determined in the
proofs below. In order to have estimates on $B_1$ and $B_2$, we
set $\bar{\psi}(X,Y) =e^{-\tilde KX_{n}}(\tilde
\psi_\delta(X,Y))^{\frac{\alpha}{2}}$, with the correspondence
given above between $X,Y$ and $x,y$.

By the regularity properties of $\psi_\delta, \phi_0$ given in
Lemma~\ref{lemmaprelbis}, it is tedious (but easy) to check that
all the terms in $D^2\Theta_{\delta}(x,y)$ are bounded except
perhaps the ones from $D^2 \bar{\psi}(X,Y)$.   Therefore, the
inequality \rec{matrixbis} can be rewritten as~: for all
$\xi,\zeta \in\r^n$
\begin{eqnarray}\label{matrix2bis}
\langle B_1\xi,\xi\rangle-\langle B_2\zeta,\zeta\rangle &\leq &
(1+O(\rho))\left[\langle P_{1} (\xi-\zeta),(\xi-\zeta)\rangle +
\langle P_{2} (\xi+\zeta),(\xi-\zeta)\rangle\right.\nonumber \\ &&
\left.+
   \langle P_{3} (\xi+\zeta),(\xi+\zeta)\rangle +
K(|\xi|^2+|\zeta|^2)\right],
\end{eqnarray}
for some constant $K$ and where
$$P_{1}= D_{YY}^2 \bar{\psi}(X,Y)\; ,\;P_{2}= D_{XY}^2 \bar{\psi}(X,Y)\;
,\;P_{3}= D_{XX}^2 \bar{\psi}(X,Y)\; .$$

Moreover, as we remark above, we can assume that $Y=x-y$ remains
bounded away from $0$ and this implies (after again tedious
computations) that $P_{3}$ is bounded as well.

Choosing
$\xi=\zeta$ in the above inequality, we first deduce that $$B_1-B_2
\leq \Ktd Id$$ with
$$\Ktd=K\left[C\left(\alpha\eta^{-1}|Y|^{\alpha
-1}(\delta^2+|Y|^2)^{1/2}+\delta|Y|^{\alpha-4}\right)+1\right]\,.$$
We next choose in \rec{matrix2bis}, $\xi=-\zeta=\hat Y.$ By using
the properties on the first and second derivatives of
$\tilde\psi_\delta$ and $\phi_0$ proved in Section 3.2, we get,
for some $K>0$
\begin{eqnarray*}
\langle (B_1-B_2){\hat Y} ,{\hat Y} \rangle  & \leq & K
Ce^{-2\tilde
KX_n}\frac{\alpha}{2}\frac{(\tilde\psi_\delta)^{\frac{\alpha}{2}
-1}}{|Y|^2}\left[\langle
D^2_{YY}(\tilde\psi_\delta)Y,Y\rangle+(\frac{\alpha}{2}
-1)\frac{\langle D_Y\tilde\psi_\delta,
Y\rangle^2}{\tilde\psi_\delta}\right]+K  \nonumber\\
   & \le &
-KC{\alpha}(1-\alpha)\frac{(\tilde\psi_\delta)^{\frac{\alpha}{2}}}{|Y|^
2}+ C\alpha O(\eta)|Y|^{\alpha-2}
   + O(\delta)|Y|^{\alpha-6} + K  \label{matrix3}
\,.\end{eqnarray*}
If $(e_i)_{1\le i\le n-1}$ are $(n-1)$ vectors such that
$(e_1,\ldots,e_{n-1},\hat Y)$ is an orthonormal basis of
$\R,$ we know that
$$ \mbox{Tr}~(B_1-B_2) = \sum_{i=1}^{n-1} \langle (B_1-B_2)e_i
,e_i\rangle +  \langle (B_1-B_2){\hat Y},{\hat Y}\rangle\, ,$$
and by combining the above estimates, we deduce
\begin{equation}\label{estitrace}
\mbox{Tr}~(B_1-B_2)\le -CK\alpha(1-\alpha)
\frac{(\tilde\psi_\delta)^{\frac{\alpha}{2}}}{|Y|^2} +
O(\delta)|Y|^{\alpha-6}  +C\alpha O(\eta)|Y|^{\alpha-2}+ \Ktd\; .
\end{equation}
Finally,   for $|Y|$ small enough (i.e. for $C$ large enough), we
get
\begin{equation}\label{final}
   \mbox{Tr}~(B_1-B_2)\le
   -CK\alpha(1-\alpha)
|Y|^{\alpha-2}+ O(\delta)|Y|^{\alpha-6}+C\alpha
O(\eta)|Y|^{\alpha-2}+ \Ktd\; .
\end{equation}
Now by using the estimates on the first and second derivatives on
$\Theta_\delta$ shown in the Subsection 3.2, we get, for some
   $K>0$,
   \begin{eqnarray*}
   |p|,|q|&\ge& C\alpha\left(
   K^{-1} - K\eta \right) |Y|^{\alpha-1}+ O(\delta)|Y|^{\alpha-3}-K\, ,\\
    |p|,|q|&\le& C\alpha K|Y|^{\alpha-1}+ O(\delta)|Y|^{\alpha-1}+K \,
,\\
    |p-q|&\le& C\alpha K|Y|^{\alpha}+O(\delta)+O(|x-x_0|^3)\, ,\\
   ||B_1||,||B_2||&\le& K (1+\displaystyle\frac{1}{O(\rho)})(1
+C\alpha|Y|^{{\alpha}-2}+O(\delta))\,. \end{eqnarray*}
As in the proof of Lemma~\ref{keylemma}, we notice that we are going to
let first $\delta$ tends to $0$ for fixed $C$ and we recall that, since
we assume that $M_{L,C}>0$, $Y$ does not converge to $0$  when $\delta$
tends to $0$ for fixed $C$. The first consequence of this fact is again
that the term $O(\delta)|Y|^{\alpha-3}$ is playing no role in the lower
estimate of $|p|, |q|$ since we can choose $\delta$ as small as
necessary. The new point in the above estimate is the $\eta$-term : we
choose it sufficiently small in order to have, say, $K^{-1} - K\eta
\geq \frac{1}{2}K^{-1}$. With this choice, by the above estimates, we
have $|p|, |q| \to + \infty$ as $C \to +\infty$

We are going to use {\bf (H2)-(H3b)} with
$R=||u||_{L^\infty(B_{\overline O}(x_{0},r)}\,.$ We drop the
dependence in $R$ in the coefficients and modulus which appear in
these assumptions.  We subtract the two inequalities \rec{ineq2}
and write the difference in the following way
\begin{eqnarray}\label{diffbis}
    && F(x,u(x),p,B_1)-
F(x,u(x),p, B_2+\Ktd Id) \\ &&~~~~~~~~\leq F(y,u(y),q,B_2)
   - F(x,u(x),p, B_2+\Ktd Id)\; ,\nonumber
\end{eqnarray} and, using the fact that $B_1-B_2 \leq
\Ktd Id$, we apply {\bf (H3b)} to the left-hand side
and {\bf (H2)} to the right-hand side of (\ref{diffbis}). This
yields
   \begin{eqnarray}\label{finalbis}
   \lambda\mbox{Tr}(B_2-B_1+\Ktd Id) &  \le &
     \omega_1 (|x-y|(1+|p|+|q|) +\varpi(|p|\wedge|q|)|p-q|) ||B||
\nonumber
\\&+& \omega_{2}(\Ktd) + C_1\nonumber
\\ &+& C_2 (|p|^2+|q|^2) +C_3|x-y|(|p|^3+|q|^3) \,.   \end{eqnarray}
The estimates on the two sides of \rec{finalbis} are done in the same
way as in the proof of Theorem~\ref{theorregul}. The only difference is
a term of the form $C\alpha
O(\eta)|Y|^{\alpha-2}$ in the left-hand side of the estimate.

Taking in account this additional term, we are lead to an analogous
estimate to (\ref{finalestim}) with a right hand side of the form
$\lambda K^{-1} + O(\eta)$ instead of $\lambda K^{-1} $. We conclude in
the same way by choosing first $\eta$ small enough.\hfill\sn

\par
\medskip
\noindent{\bf Step 2 : The case when $G$ is independent of $u$ but with
a general dependence in $p$}\newline
As for the treatment of the dependence in $u$, we are going to
introduce a new variable. More precisely we introduce the function $v :
\overline O \times \r \to \r$ defined by
$$ v(x,y) := u(x) - y\; .$$

This new function is formally a solution of
$$ -\frac{\partial^2 v}{\partial y^2}+F(x, v+y , D_x v , D_{xx}^2v) =
0\quad\mbox{in $O\times \r\,$} $$
and
$$
-D_y v G\left(x, -\frac{D_x v}{D_y v}\right) =0\quad\mbox{in
$\partial O\times \r\,.$} $$
In fact, in order to justify this, one has just to be a little bit more
precise about the definition of the boundary condition. We set, for
$x\in \partial O$, $p_x \in \R$ and $p_y \in \r$
$$ \tilde G(x,p_x,p_y):= \left\{
\begin{array}{ll}
-p_y G\left(x,\displaystyle -\frac{p_x}{p_y}\right) & \mbox{if
$p_y<0$,}\\
G_\infty(x, p_x) & \mbox{if  $p_y \geq 0$.}
\end{array}\right.$$
With this notation, the boundary condition for $v$ becomes $\tilde
G(x,D_x v,D_y v)= 0$ and, because in particular of the assumptions {\bf
(G3)-(G4)}, it is rather easy to show that $\tilde G$, in addition to
be homogeneous of degree 1 in $(p_x,p_y)$, satifies {\bf (G1)-(G2)}.

On the other hand, the assumptions on the equation can be checked
easily and therefore the conclusion follows from Step 1.

\smallskip
\noindent{\bf Step 3 : The general case}\newline
In order to treat the dependence in $u$, as mentioned above, we use the
method of the second step of the proof of Theorem~\ref{theorregul}. We
are not going to give all the details since they are essentially the
same. We just want to point out that in order to take care of the
dependence in $u$ and to be sure that the transformed boundary
condition actually satisfies {\bf (G3)-(G4)}, one has first to
introduce the function $G_R$ defined by $R>0$ large enough and for
$x\in \partial O$, $u\in \r$ and $p\in \R$, by
$$G_R(x,u,p):=\left\{
\begin{array}{ll}
G(x,u,p) & \mbox{if  $|u|\leq R$,}\\
G(x,-R,p) & \mbox{if  $u\leq -R$,}\\
G(x,R,p) & \mbox{if  $u \geq R$.}\end{array}\right.$$
Clearly, if $R$ is large enough, $u$ is still a solution of the Neumann
problem with $G_R$ and this transformation prevents difficulties with
the behavior of $G$ in $u$ for $|u|$ large.

And the proof of Theorem~\ref{theorregulbis} is complete.\hfill\sn

\medskip
\begin{remark} We want to point out that, in {\bf (H2)}, the $\varpi$
term is needed only because we want to obtain local estimates~: this is
clear in the proof since, in (\ref{finalestim}), this $\varpi$ term is
used to take care of $Q_2$ which comes from the localization term.
Therefore, in the case of global estimates, the same result holds with
$\varpi \equiv 1$. This remark can be used either on bounded domains or
in unbounded domains where, if the equation and the boundary condition
satisfy suitable uniform properties, $L$ can be taken as small as we
want (a ``mild'' localization) and the same effect occurs.
\end{remark}

\section{The construction of the test-functions}\label{CTF}

In this section,  we provide the proof of Lemma \ref{lemmaprel} and
\ref{lemmaprelbis}. In particular we show how we construct the
functions $\Theta_0$ and $\Theta_{\delta}$ which are used in the
proof of Lemma~\ref{keylemma} and \ref{flatbound}. We will
consider separately as in Section 1 the case of linear and
nonlinear boundary conditions.

\subsection{The test-function for linear boundary conditions}

In this subsection we consider the case
$$\frac{\partial u}{\partial \gamma} + g(x) = 0 \quad\hbox{on  }\partial
O\; ,$$ where $\gamma\colon\partial O\to\r^n$ is  a locally
Lipschitz continuous vector field such that $\langle\gamma(x),
n(x)\rangle\ge \nu>0$ for any $x\in \partial O$, and $g\colon\partial
O\to \R $ is either a
locally Lipschitz  continuous or a locally H\"older continuous
scalar function.

According to assumption {\bf (H3a)}, there exists a function
$A(\cdot)\in C_{loc}^{0,1}(\partial O,{\cal{S}}^n)$ such that,
for any $x \in \partial O$, $n(x)=A(x)\gamma(x)$ and $A(x)\ge c_0 Id$,
for some constant $c_0>0$ and such that (\ref{ndpsod}) holds. Of
course, this last property is the most important information in {\bf
(H3a)}, the existence of such $A$ without the connection with the
ellipticity of the equation, being easy to show.

We may assume  without loss of generality, that $\langle\gamma(x),
n(x)\rangle=1$ for any $x \in \partial O$, otherwise we change $\gamma$
in
$\frac{\gamma(x)}{\langle\gamma(x), n(x)\rangle}$, $g$ in
$\frac{g(x)}{\langle\gamma(x), n(x)\rangle}$ and $A(x)$ in
$\langle\gamma(x), n(x)\rangle A(x)$; these transformations do not
change the properties of $\gamma$ and $g$.

\medskip
{\noindent\bf Proof of Lemma \ref{lemmaprel}.}
As in the proof of the comparison result for this kind of problems (cf.
Barles \cite{Ba1}), we are going to use regularizations of  $A$
and  $g$. To do so, it is convenient to introduce the following
lemma whose proof is classical and therefore left to the reader.

   \begin{lemma}\label{basicreg} Assume that $\rho\in{{D}}(\R)$, $\rho\ge
0$, supp$(\rho) \subset B(0,1)$ and $\int_{\R} \rho(y)\,dy =1.$ If
$f \in C^{0,\beta}(\R) $ for some $0< \beta \le 1$, and $f$ is
bounded, then the function $\tilde f : \R \times [0,+\infty) \to
\r$ defined, by
   $$ \tilde f(x,\eps) := \int_{\R}
f(z)\rho(\frac{x-z}{\e})\frac{1}{\e^n}\,dz\quad\hbox{for }x\in \R,
\;\eps > 0 ,$$

$$ \tilde f(x,0) = f(x)\quad\hbox{for }x\in \R,$$
   is   in $ C^{0,\beta}(\R\times [0,+\infty))$.
   Moreover, the function $\tilde f$ is $C^2$ in $\R\times (0,+\infty)$
with

$$|D_x \tilde f(x,\eps)|, |D_\eps \tilde f(x,\eps)|\le K\eps^{\beta-1}$$
$$ |D^2_{xx}  \tilde f(x,\eps)|, |D^2_{x\eps}  \tilde f(x,\eps)|,
|D^2_{\eps\eps}
   \tilde f(x,\eps)|\leq K\eps^{\beta-2}\quad \hbox{in  }\R\times
   (0,+\infty)$$
    for some constant $K$ depending only on $\rho$, the $L^\infty$ and
the H\"older
norm  of $f$.
\end{lemma}

   \medskip

   {\noindent\bf Step 1.} The functions $A$  and  $g$
   and their regularizations.

    We first extend
$g$ and $A$    to $\R$; we still denote by $g$ and $A$ these
extensions. We may assume that these extensions are respectively
in in $C^{0,\beta}(\R)$ and $C^{0,1}(\R)$. For some function
$\rho$ satisfying the properties of Lemma~\ref{basicreg} (which is
chosen and fixed from now on), we consider the functions $\tilde
A$ and $\tilde g$ associated to $A$ and $g$ as in this lemma.
Finally, we introduce, for some   $\delta \geq 0$, the following
quantity which is defined for $\xi  \in \R$ by
   $$
\Lambda_\delta(\xi)=\left(\delta^2+\vert \xi\vert^2\right)^{1/2}\, ,
$$
and we set
$$
   \tilde A_\delta(x,\xi)= \tilde A (x,\Lambda_\delta(\xi)) \quad , \quad
   \tilde g_{\delta} (x, \xi)= \tilde g(x, \Lambda_\delta(\xi))\; .
$$
According to Lemma~\ref{basicreg}, these functions are $C^2$ as
long as $\delta >0$.

We also observe that
$$ |D_\xi\Lambda_\delta(\xi)|\le 1\;,\;
|D^2_{\xi\xi}\Lambda_\delta(\xi)|\le
K\Lambda_\delta^{-1}\,.
$$

  \medskip

{\bf\noindent Step 2.} Construction of the functions $\psi_0$,
$\psi_\delta$ and their main properties.

\smallskip
For $\delta \geq 0$, we introduce the following function, for
$X,Y\in \R$ and $T>0$
\begin{equation}\label{psidelta}
\tilde \psi_\delta (X,Y,T)= \langle \tilde
A_\delta(X,Y)Y,Y\rangle+ K_1\delta (2M - T)\; ,
\end{equation}
where $ K_1>0$ is a constant  to be chosen later and $M$ is chosen
so that $2M - T$ remains bounded. Moreover we set, for $x$ and $y$
in a suitable neighborhood of $x_0$ $$ \psi_\delta (x,y) := \tilde
\psi_\delta \left( \frac{x+y}{2},x-y,d(x)+d(y) \right)\; .$$

We observe that, as $\delta\to 0$, $\tilde \psi_\delta$ and
$\psi_\delta$ converge locally uniformly respectively
   to $\tilde \psi_0$ and $\psi_0$.
   Depending on the simplicity, we provide below result either on $\tilde
\psi_\delta$ or $\psi_\delta$,
   the translation from one to the other being straightforward. Most of
the time we will use $\tilde \psi_\delta$.

In the sequel $K>0$ denotes a constant which may
vary from line to line but which depends only on the data of the
problem and is independent of the small parameter $\delta.$

\begin{prop}
We have, for any $X, Y \in \R$
\begin{eqnarray}
& &K^{-1}|Y|^2\le \tilde \psi_0(X,Y)\le K |Y|^2\,,\label{psiopos}\\
& &  K^{-1}|Y|^2\le \tilde \psi_\delta (X,Y)\le K
|Y|^2+K\delta\,.\label{psideltapos}
\end{eqnarray}
\end{prop}

The proposition  is straightforward consequence of the fact that
$A(x)\ge c_0 Id$ for all $x\in\r^n$. Next we examine the regularity
properties and the estimates on
$\tilde \psi_\delta$ and $\psi_\delta$.

\begin{prop} \label{estimatetildepsi}
We have, for any $X,Y\in \R$, $T\in \r$
\begin{eqnarray*}
&& |D_{Y}\tilde \psi_\delta (X ,Y,T)|\le K|Y|\;,\; |D_X \tilde
\psi_\delta(X ,Y,T)|\le K|Y|^2\;, \\ && | D^2_{{Y}{Y}}
\tilde\psi_\delta (X,Y,T)|\le K\;,\;|D_{{X}{X}}^2 \tilde
\psi_\delta (X,Y,T)|\leq K|Y|\,,
  \\ &&
|D_{{X}{Y}}^2 \tilde \psi_\delta (X ,Y,T)|\leq K |Y|\,, \nonumber \\
&&  |D_T\tilde \psi_\delta (X ,Y,T)|\leq K\delta\, , \,D_{TT}^2 \tilde
\psi_\delta=D_{TX}^2 \tilde
\psi_\delta=D_{TY}^2 \tilde \psi_\delta=0\;.
\end{eqnarray*}
Moreover
\begin{eqnarray*}
\langle D_Y  \tilde\psi_\delta (X,Y,T), Y\rangle&=&
2\tilde\psi_\delta+O(|Y|^3)+O(\delta)\, ,\\
\langle D^2_{YY}  \tilde\psi_\delta Y,Y\rangle&=&
   2 \tilde\psi_\delta+O(|Y|^3)+O(\delta)\,,
\end{eqnarray*}
as $Y\to 0$ and $\delta\to 0\,.$
\end{prop}

The proof of these estimates is tedious but straightforward: the
main reason to provide Lemma~\ref{basicreg} and to write $\tilde
A_\delta$ with a dependence in $x$ and $\Lambda_\delta(\xi)$ was
to have a simple way to check these computations.
\par
{\noindent \bf Proof of Proposition~\ref{estimatetildepsi}.}  We have
\begin{eqnarray*}
    D_X \tilde\psi_\delta (X,Y,T) &=&\langle D_X\tilde A_\delta
    Y,Y\rangle\, ,\\
D_Y \tilde\psi_\delta (X,Y,T) &=& \langle D_Y\tilde A_\delta
    Y,Y\rangle+2\tilde A_\delta Y\, ,\\
D^2_{XX}\tilde\psi_\delta (X,Y,T) &=&\langle D_{XX}\tilde A_\delta
    Y,Y\rangle \,,
\\
D^2_{YY}\tilde\psi_\delta (X,Y,T) &=&\langle D_{YY}\tilde A_\delta
    Y,Y\rangle+2\tilde A_\delta+2D_Y\tilde A_\delta Y\,.
\end{eqnarray*}
We premise some useful estimates on the first and second
derivatives of $\tilde A_\delta$. By using Lemma \ref{basicreg}
and the estimates on the first and second derivatives of
$\Lambda_\delta$, we have
\begin{eqnarray*}\label{proptildeA}
& & D_X\tilde A_\delta (X,Y)=D_{\Lambda_\delta}\tilde A =O_Y(1)\, ,
\\&& D^2_{XX}\tilde A_\delta (X,Y)= D^2_{\Lambda_\delta}\tilde A
=\Lambda_\delta^{-1}O_Y(1) \, , \nonumber
\\
& & D_Y\tilde A_\delta (X,Y)=D_{\Lambda_\delta}\tilde
A\frac{\partial\Lambda_\delta}{\partial Y}
= O_Y(1)\, , \nonumber\\
& & D^2_{YY}\tilde A_\delta (X,Y)  = (D_{\Lambda_\delta}\tilde
A)\frac{\partial^2\Lambda_\delta}{\partial^2 Y}
+D^2_{\Lambda_\delta}\tilde
A\left(\frac{\partial\Lambda_\delta}{\partial Y}\right)^2=
\Lambda_\delta^{-1}O_Y(1)\,. \nonumber
\end{eqnarray*}

Using the estimates on the first and second derivatives of properties
of $\tilde A_\delta$,  we obtain easily all the estimates on the
derivatives of $\tilde\psi_\delta $ and the second part of
Proposition~\ref{estimatetildepsi}.\hfill\sn

   \medskip
We turn to the properties of $\psi_\delta$ with respect to the boundary
condition.
\begin{prop}\label{bctildepsi}
If $|x-y|$ is small enough and $ K_1$ is large enough, then we
have, for some $K>0$
\begin{eqnarray}
&& \langle D_x\psi_\delta(x,y),\gamma (x)\rangle >-K|x-y|^2~~~\mbox{if
$x\in\partial O\,,$}\label{bcx2}
\\
&& \langle-D_y\psi_\delta(x,y),\gamma(y)\rangle<K|x-y|^2~~~\mbox{if
$y\in\partial O\,.$}\label{bcy2}
\end{eqnarray}
\end{prop}
{\noindent\bf Proof of Proposition \ref{bctildepsi} .}
We only check \rec{bcx2} the other case being
similar.

By a direct computation, we have
\begin{eqnarray*}
D_x\psi_\delta(x,y)&=& 2\tilde A_\delta
(x-y)+\langle(\frac{D_X\tilde A_\delta}{2}+D_Y \tilde A_\delta)
Y,Y\rangle-K_1\delta Dd(x)\,\\ &=& 2\tilde A_\delta (x-y)+O(\vert
x-y\vert^2) +K_1\delta n(x)
\end{eqnarray*}
But we recall that
   $$||\tilde A_\delta-A||\le K(|x-y|+\delta)$$
and since $x\in\partial O$, by the regularity of the boundary, we have
$$
\langle n(x),(x-y)\rangle=d(y)+O(|x-y|^2)~~~\mbox{as $|x-y|\to 0\;.$}
$$
Thus    if $K_1$ is
large enough  we have
\begin{eqnarray*}\label{dxpsi}
\langle D_x\psi_\delta(x,y),\gamma(x)\rangle&\ge & \langle 2
A(x-y),\gamma\rangle -K|x-y|(|x-y|+\delta) -K|x-y|^2 +K_1\delta\, , \\
&\ge&\langle 2(x-y), n(x)\rangle -K\delta -K|x-y|^2+K_1\delta \, , \\
&\ge& -K|x-y|^2\, .
\end{eqnarray*}

   \hfill\sn

\medskip
{\bf\noindent Step 3.} Construction of the functions $\chi_0$ and
$\chi_\delta.$

In the same way as above, we set for $\delta \geq 0$
\begin{equation}\label{chidelta}
\tilde \chi_\delta(X,Y,T,Z)=\tilde g_\delta (X,Y)Z +
K_2\delta^\beta (2M-T)\; ,
\end{equation}
where $ K_2>0$ is a constant to be chosen later and $M>0$ is chosen as
above. We also set
$$ \chi_\delta (x,y) := \tilde \chi_\delta \left(
\frac{x+y}{2},x-y,d(x)+d(y),d(x)-d(y) \right)\; .$$

One can easily check that $\tilde \chi_\delta$ and $ \chi_\delta$
converges locally uniformly respectively
   to $\tilde \chi_0$ and $ \chi_0$ as $\delta \to 0$.

As in the previous step, we first consider the regularity properties of
$\tilde \chi_\delta$.

\begin{prop}\label{estimatetildechi}
For every $\delta>0$ we have
$$ |D_X \tilde \chi_\delta (X,Y,T,Z) | \leq  K
\Lambda^{\beta-1}_\delta |Z| \; , \; |D_Y \tilde \chi_\delta
(X,Y,T,Z) | \leq K \Lambda^{\beta-1}_\delta |Z|\; ,$$
$$|D^2 \tilde \chi_\delta (X,Y,T,Z)|\leq  K
\Lambda^{\beta-2}_\delta |Z| \; .$$
\end{prop}
{\noindent\bf Proof of Proposition \ref{estimatetildechi}.}
Again the tedious computations are simplified by the way we write down
$\tilde g_\delta$.
We first observe that
\begin{eqnarray*}
D_X\tilde g_\delta&=&D_X\tilde g=\Lambda^{\beta-1}_\delta O_Y(1)\, , \\
D^2_{XX}\tilde g_\delta&=&D_{XX}\tilde g
= \Lambda^{\beta-2}_\delta O_Y(1)\, ,\\
D_Y\tilde g_\delta&= &D_{\Lambda_\delta}\tilde g
\frac{\partial\Lambda_\delta}{\partial Y}=\Lambda^{\beta-1}_\delta
O_Y(1)\, ,\\
D^2_{YY}\tilde g_\delta&=& D^2_{\Lambda_\delta}\tilde g
\left(\frac{\partial\Lambda_\delta}{\partial
Y}\right)^2+D_{\Lambda_\delta}\tilde g
\frac{\partial^2\Lambda_\delta}{\partial^2 Y}
=\Lambda^{\beta-2}_\delta O_Y(1)\,.
\end{eqnarray*}
Then, by the definition of $\tilde\chi_\delta$, we deduce
\begin{eqnarray*}
D_X\tilde\chi_\delta(X,Y,T,Z)&=& D_X\tilde g_\delta Z\le K
\Lambda^{\beta-1}_\delta |Z|\, ,
    \nonumber\\
D_Y\tilde\chi_\delta(X,Y,T,Z)&=&D_Y\tilde g_\delta Z\le K
\Lambda^{\beta-1}_\delta |Z|\, , \\
D^2_{XX}\tilde\chi_\delta(X,Y,T,Z)&=& D^2_{XX}\tilde g_\delta Z
\le K
\Lambda^{\beta-2}_\delta |Z| \, ,\nonumber\\
D^2_{YY}\tilde\chi_\delta(X,Y,T,Z)&=&D^2_{YY}\tilde g_\delta
Z\le K \Lambda^{\beta-2}_\delta |Z|\,.\end{eqnarray*}
~\hfill\sn

Then we consider the boundary condition.
\begin{prop}\label{bctildechi}
If $|x-y|$ is small enough and $K_2$ is large enough, then we have, for
some $K>0$
\begin{eqnarray}
&& \langle D_x\chi_\delta(x,y),\gamma (x)\rangle  + g(x) >  -
K|x-y|^{\beta}+K\delta^{\beta}~~\mbox{if
$x\in\partial O\, ,$}\label{bcchix2}
\\
&& \langle -D_y\chi_\delta(x,y),\gamma(y) \rangle+ g(y)<
K|x-y|^{\beta}-K\delta^{\beta}~~\mbox{if
$y\in\partial O\, .$}\label{bcchiy2}
\end{eqnarray}
\end{prop}

{\noindent\bf Proof of Proposition \ref{bctildechi}
.}
Again we only check the first property \rec{bcchix2}. We first notice
that, by the definition and properties of $\tilde
g_\delta$,
$$
|\tilde g_\delta(\frac{x+y}{2},x-y)- g(x)|\le
K(|x-y|^\beta+\delta^\beta).
$$
On the other hand, we have
$$
D_x \chi_\delta(x,y)=
   \tilde g_\delta Dd(x) +(\frac{1}{2} D_X \tilde g_\delta+D_Y \tilde
g_\delta ) (d(x)-d(y))
-  K_2\delta^\beta Dd(x).
$$
Therefore since $|x-y|$ is small,  we obtain for $K_2$ large
enough
\begin{eqnarray*}
\langle D_x\chi_\delta(x,y),\gamma(x)\rangle+g(x) & \ge & [g(x) -\tilde
g_\delta (\frac{x+y}{2},x-y)] \\&+& \langle
( \frac{1}{2} D_X \tilde g_\delta+D_Y \tilde g_\delta ),\gamma
(x)\rangle (d(x)-d(y)) \\
    &+& K_2 \delta^\beta\\
&\ge&  -K|x-y|^\beta+K\delta^\beta \; ,
\end{eqnarray*}
Thus we have shown \rec{bcchix2}.~\hfill\sn

\medskip
Now we can prove the following result.
\begin{prop}\label{bctheta}
If $|x-y|<<1$, then  we have
\begin{eqnarray}
&& \langle D_x\Theta_{\delta}(x,y), \gamma(x)\rangle+ g(x)>0~~\mbox{if
$x\in\partial
O\, ,$}\label{hbcxn3}
\\
&& \langle-D_y\Theta_{\delta}(x,y), \gamma(y)\rangle+g(y)<0~~\mbox{if
$y\in\partial
O\,.$}\label{hbcyn3}
\end{eqnarray}
   \end{prop}
{\noindent\bf Proof of Proposition \ref{bctheta}.}
Again we only check the first property \rec{hbcxn3}.   We first
observe that, because of the properties of $\psi_\delta$ and since
$d(x)= 0$, we have
\begin{eqnarray*}
\langle D_x(e^{-\tilde K(d(x)+d(y))}
(\psi_\delta(x,y))^{\frac{\alpha}{2}}),\gamma(x)\rangle
&=&\tilde K e^{-\tilde Kd(y)} \tilde \psi_{\delta}(x,y) \\
&+&\frac{ \alpha}{2} e^{-\tilde Kd(y)}
(\psi_\delta(x,y))^{\frac{\alpha}{2}-1})\langle D_x
\psi_\delta(x,y),\gamma (x)\rangle \\&\ge&
\tilde K e^{-\tilde Kd(y)} |x-y|^\alpha - K   e^{-\tilde
Kd(y)}|x-y|^\alpha.
\end{eqnarray*}
Thus, if $\tilde K$ is sufficiently large,  we obtain
$$\langle D_x \left[ e^{-\tilde  K(d(x)+d(y))}
(\psi_\delta(x,y))^{\frac{\alpha}{2}}\right],\gamma(x)\rangle> \bar K
|x-y|^\alpha\; ,$$ for some constant $\bar K$. Similarly one can shows
that
$$\langle D_x \left[ e^{-\tilde  K(d(x)+d(x_0))}
(\psi_\delta(x,x_0)) \right],\gamma(x)\rangle> \bar K |x-x_0|^2\,.
$$

  Then if $|x-y|<<1$ and $C>0$ is large enough, we have
\begin{eqnarray*}
\langle D_x\Theta_{\delta}(x,y),\gamma(x)\rangle+g(x)&=&\langle CD_x
[e^{-\tilde
K(d(x)+d(y))} (\psi_\delta(x,y))^{\frac{\alpha}{2}})],\gamma(x)\rangle
\\&+&\langle LD_x [e^{-\tilde
Kd(x)}\psi_\delta(x,x_0)],\gamma(x)\rangle\\
&+&\langle D_x\chi_\delta(x,y),\gamma\rangle+ g(x)\\
&\ge& C K
e^{-\tilde{K}(d(x)+d(y))}|x-y|^\alpha-K|x-
y|^\beta+K\delta^{\beta}>0.~~~~~\mbox{$\sn$}
\end{eqnarray*}

\medskip
The proof of Lemma \ref{lemmaprel} is obtained by combining
Propositions~\ref{estimatetildepsi},
\ref{estimatetildechi} and \ref{bctheta}.

\medskip
We conclude this section with the following lemma.
\begin{lemma}\label{lemmaAY2}
If $P = A^{-1}(x)D_x\Theta_\delta (x,y)$, we have
$${\hat P}={\hat Y}+o_Y(1) + O(\delta)
\quad\mbox{as $|Y|\to 0$ and $\delta \to 0$.}
$$
\end{lemma}
{\noindent \bf Proof.} We just give a sketch of proof. We first
observe that
$$||A_\delta(\frac{x+y}{2},x-y)-A(x)||\le K(|x-y|+\delta).$$
Therefore, by direct computations, we obtain
\begin{eqnarray*}
D_x\Theta_\delta(x,y)&=&Ce^{-\tilde
K(d(x)+d(y))}\left(\tilde\psi_\delta\right)^{\alpha/2
-1}\left[AY+O(|Y|^2)+O(\delta)\right]\\
&=&Ce^{-\tilde
K(d(x)+d(y))}\left(\tilde\psi_\delta\right)^{\alpha/2
-1}\left[AY+O(|Y|^2)
+ +O(\delta)\right]\,.
\end{eqnarray*}
Thus we get
\begin{eqnarray*}
\frac{A^{-1}(x)D_x\Theta_\delta(x,y)}{|A^{-1}(x)D_x\Theta_\delta(x,y)|}
&=&\frac{Y+O(|Y|^2)+O(\delta) }{|Y+O(|Y|^2)+O(\delta)|}\\
&=& {\hat Y}+ o_Y(1) + O(\delta)~~~\mbox{as $|Y|\to 0\,,\delta\to
0\,,$}
\end{eqnarray*}
which gives the desired result.~\hfill\sn

\subsection{The test-function for nonlinear boundary conditions}

We recall that we have to build this test-function in the case when the
function $G$ is independent of $u$ and
homogeneous of degree $1$ with respect to $p$.

We first extend the function $G(x,p)$ to $\R \times \R$ and we may
assume that all the
properties of $G$ are still satisfied in ${\cal{V}}\times
\r^n$ where  ${\cal{V}}$ is a neighborhood of $\partial O$. The
properties of $G$ imply that, for every $x\in
{\cal{V}},$ $p\in\r^n$ there exists a unique solution
$t=\varphi(x,p)$ of the equation
\begin{equation}
G(x,p+t n(x))=0.
\end{equation}
One can verify that $\varphi$ is still homogeneous of degree $1$
and satisfies ${\bf(G2)}.$

It is not restrictive to reduce to the case when the boundary is flat
and more
precisely $O=\{x_n>0\}$  and $  \partial O=\{x_n = 0\}.$

\medskip
{\noindent\bf Proof of Lemma \ref{lemmaprel}.}\par

{\noindent\bf Step 1.} The function $\varphi$ and its regularization.

In order to regularize the function $\varphi$, we
first extend  it  to $\r^n$ and  we still denote by $\varphi$ this
extension. We may assume that this extension satisfies ${\bf (G2)}$.

We introduce, for  $\delta, \eta > 0$, the following quantity which is
defined for $\xi  \in \R$ by
$$
\Gamma_\delta(\xi)=\eta\left(\frac{\delta^2+(\xi_n)^2}{\delta^2+|\xi|^2}
\right)^{1/2}\,,
$$
and we set, for $x, \xi \in \R$,
$$
   \tilde \varphi_\delta(x,\xi)= \tilde \varphi
   (x,{\hat \xi},\Gamma_\delta(\xi))\; ,
$$
where $\tilde \varphi$ is defined as in Lemma~\ref{basicreg} and ${\hat
\xi}= \frac{\xi}{|\xi|}$.

We first observe that the following estimates, which are used
extensively in the sequel, hold$$|D_\xi\Gamma_\delta|\le
\frac{K\eta}{(1+|\xi|^2)^{1/ 2}}\;,\; |D_{\xi\xi}\Gamma_\delta|\le
\frac{K\eta}{(1+|\xi|^2)^{1/ 2}(1+|\xi_{n}|^2)^{1/ 2}}\,. $$

\medskip

{\bf\noindent Step 2.} Construction of the functions $\psi_0$,
$\psi_\delta,$   and their main properties.

\smallskip
For $\delta \geq 0$, we introduce the following function, for
$X,Y\in \R$ with as above $Y=(Y_1, \ldots ,Y_n)$ and $X=(X_1, \ldots
,X_n)$
\begin{equation}\label{psIdelta}
\tilde \psi_\delta (X,Y)= |Y|^2-2\tilde
\varphi_\delta(X,Y)|Y|Y_{n}+2A_{1}Y_{n}^2+
K_1\delta(R-X_n)\;,
\end{equation}
with $A_1, K_1>0$  constants to be chosen later. The constant
$R>0$ has to be chosen in order that the term $R - X_n$ remains
positive; this does not create any problem since we argue locally.
Moreover we set, for $x$ and $y$ in a  suitable neighborhood of
$x_0$
$$ \psi_\delta (x,y) := \tilde \psi_\delta \left( \frac{x+y}{2},x-y
\right)\; .$$

We observe that, as $\delta\to 0$, $\tilde \psi_\delta$ and
$\psi_\delta$ converge locally uniformly respectively to $\tilde
\psi_0$ and $\psi_0$. As in previous subsection   we provide below
result either on $\tilde \psi_\delta$ or $\psi_\delta$, the
translation from one to the other being straightforward. Most of
the time we will use $\tilde \psi_\delta$.

In the sequel $K>0$ will denote a nonnegative constant which may
vary from line to line but which depends only on the data of the
problem and is independent of the small parameters $\delta$ and
$\eta$

\begin{prop}
If  $A_1>0$ is large enough, we have, for any $X,Y \in \R$
\begin{eqnarray}
& &K^{-1}|Y|^2\le \tilde \psi_0(X,Y)\le K
|Y|^2\,,\label{psioposbis}\\ & &  K^{-1}|Y|^2\le \tilde
\psi_\delta (X,Y)\le K |Y|^2+K\delta\,.\label{psideltaposbis}
\end{eqnarray}
\end{prop}

We skip the proof of this proposition which is straightforward~:
it is based only on Cauchy-Schwarz's inequality to control the
$\tilde \varphi_\delta$ or the $\tilde \varphi_0$ term and on the
fact
that $\tilde \varphi_\delta$ and $\tilde \varphi_0$ are bounded.

Next we examine the regularity properties and the estimates on
$\tilde \psi_\delta$ and $\psi_\delta$.

\begin{prop} \label{estimatetildepsibis} If the constant $\eta>0$ is
chosen small enough and $A_1$, $K_1$ large enough, then, for any
$\delta$ small enough and for all $X,Y \in \R$, we have
$$ |D_X \tilde \psi_\delta (X,Y) | \leq K|Y|^2+K\delta\; ,
   \; |D_Y \tilde \psi_\delta (X,Y) | \leq K|Y|\; ,$$
   $$ \langle D_{Y}\tilde \psi_\delta (X,Y), Y\rangle =
2\tilde\psi_\delta(X,Y)+O(\eta) O(|Y|^2) +O(\delta)\; ,$$
$$  \langle D^2_{YY}\tilde \psi_\delta (X,Y) Y,Y\rangle =
2\tilde\psi_\delta(X,Y)+O(\eta) O(|Y|^2) +O(\delta)\;,$$
as $Y\to 0,$ $\eta\to 0$ and  $\delta\to 0\,,$
   $$|D_{XX}^2 \tilde \psi_\delta (X,Y)|\leq \frac{K}{\eta}|Y|(\delta^2
+|Y|^2)^{1/2}\; , \;
|D_{XY}^2 \tilde \psi_\delta (X,Y)|\leq K |Y|  \; , \; |D_{YY}^2
\tilde \psi_\delta (X,Y)|\leq K.$$
\end{prop}
{\noindent\bf Proof.} The proof of these estimates is tedious but
straightforward.  Lemma~\ref{basicreg} and the way we write
$\tilde \varphi$ with a dependence in $x$, $\hat \xi$ and
$\Gamma_\delta(\xi)$ is
a simple way to check these computations.

By direct computations, we have
\begin{eqnarray*}
    D_X \tilde\psi_\delta (X,Y)&=& -2D_X\tilde \varphi_\delta
|Y|Y_n-K\delta
e_n \, ,\\
D_Y \tilde\psi_\delta (X,Y) &=& 2Y-D_Y\tilde \varphi_\delta |Y|Y_n \\ &&
-2\tilde \varphi_\delta\frac{Y}{ |Y|}Y_n-2\tilde \varphi_\delta
|Y|e_n+2A_1Y_ne_n\, ,\\ D^2_{XX}\tilde\psi_\delta (X,Y) &=&
-2D_{XX}\tilde \varphi_\delta |Y|Y_n \, ,
\\
D^2_{YY}\tilde\psi_\delta (X,Y) &=& 2Id-2D_{YY}\tilde \varphi_\delta
|Y|Y_n \\
&& -4D_Y \tilde \varphi_\delta\frac{Y}{ |Y|}Y_n-4D_Y \tilde
\varphi_\delta {
|Y|}e_n\, ,\\
&& -4 \tilde \varphi_\delta\frac{Y}{ |Y|}e_n-4 \tilde
\varphi_\delta D_Y\frac{Y}{ |Y|}Y_n+2A_1Y_n\; ,
\end{eqnarray*}
where $e_n$ is the $n$-th vector the canonical basis of $\r^n$. In
this case, we also have $n(x)=-e_n,$ for all $x\in\partial O.$

We estimate  the first and second derivatives of $  \tilde
\varphi_\delta$. By using Lemma \ref{basicreg} and the estimates on the
first and second derivatives of $\Gamma_\delta$, we have
\begin{eqnarray*}
  |D_X\tilde  \varphi_\delta (X,Y)| & = & |D_X\tilde
\varphi(X,\frac{Y}{|Y|},\Gamma_\delta(Y))|\le K \, ,
\\
  | D^2_{XX}\tilde \varphi_\delta (X,Y)| & = &|D^2_{XX}\tilde
\varphi(X,\frac{Y}{|Y|},\Gamma_\delta(Y))|\le \frac{K}{\Gamma_\delta}
\, ,
\\
|D_Y\tilde \varphi_\delta (X,Y) | & = &|D_{\xi}\tilde \varphi
D_Y\frac{Y}{|Y|}+D_{\zeta}\tilde  \varphi D_Y\Gamma_\delta |\le
\frac{K}{|Y|} \, , \\
  |D^2_{YY}\tilde \varphi_\delta
(X,Y) & = & |D^2_{\xi\xi}\tilde  \varphi\left(D_Y\frac{Y}{|Y|}\right)^2
+D_{\xi}\tilde  \varphi D^2_{YY}\frac{Y}{|Y|} +2D^2_{\xi\zeta}\tilde
\varphi D_Y\frac{Y}{|Y|} \otimes
D_Y\Gamma_\delta \\ &&
+ D^2_{\zeta\zeta}\tilde
\varphi\left(D_Y\Gamma_\delta \right)^2+D_{\zeta}\tilde \varphi
D^2_{YY}\Gamma_\delta |\\&\le &
   \frac{K}{  |Y | |Y_n |}(\frac{1}{\eta}+1)\;.
\end{eqnarray*}

\smallskip
\noindent By combining the above estimates we obtain
$$
   | D_X  \tilde\psi_\delta (X,Y)| \le  K|Y|^2+O(\delta)\, , \,
|D_Y\tilde\psi_\delta (X,Y)|\le K|Y|\, , $$
$$
   | D_{XX}^2  \tilde \psi_\delta (X,Y)| \le \eta^{-1}
K|Y|(|Y|^2+\delta^2)^{1/2}
\, , \, |D_Y\tilde\psi_\delta (X,Y)|\le K|Y|\,. $$

\medskip
\noindent Next we estimate $\langle D_Y
\tilde\psi_\delta (X,Y)  , Y\rangle$ and $\langle D^2_{YY}
\tilde\psi_\delta Y,Y\rangle.$ Tedious but straightforward
computations show that
\begin{equation}\label{usful1}
|\langle D_Y\Gamma_\delta, Y\rangle|\le K\eta\,,\quad |\langle
D_{YY}\Gamma_\delta Y,Y\rangle|\le K\eta\,.
\end{equation}
Moreover
\begin{eqnarray*}
\langle D_Y  \tilde\psi_\delta (X,Y), Y\rangle&=&
   2 |Y |^2-4\tilde \varphi_\delta  |Y |Y_n+2A_1 Y^2_n\\
   &&- D_{\zeta}\tilde \varphi_\delta (\langle
D_Y\Gamma_\delta,Y\rangle) |Y |Y_n\\
   &=& 2\tilde\psi_\delta +O(\delta)+O(\eta)|Y|^2\;;
\end{eqnarray*}

\begin{eqnarray*}
\langle D^2_{YY}  \tilde\psi_\delta Y,Y\rangle&=&2 |Y |^2-4\tilde
\varphi_\delta  |Y |Y_n+2A_1 Y^2_n
\\
&&-2D_{\zeta}\tilde \varphi_\delta\langle D^2\Gamma_\delta Y,Y\rangle
|Y |Y_n- 4D_{\xi}\tilde \varphi_\delta(\langle D\Gamma_\delta,
Y\rangle)  |Y |Y_n\\
&&- D^2_{\xi\xi}\tilde \varphi_\delta(\langle D\Gamma_\delta,
Y\rangle)^2  |Y
|Y_n\\ &=& 2\tilde\psi_\delta +O(\delta)+O(\eta)|Y|^2\;.
\end{eqnarray*}

These properties complete the proof of
Proposition~\ref{estimatetildepsibis}.\hfill\sn

   \medskip
We turn to the properties of $\psi_\delta$ with respect to the
boundary condition.
\begin{prop}
If $|x-y|$ is small enough and $ K_1$ is large enough, then we
have, for some $K>0$
\begin{eqnarray}
&& G(x,D_x\psi_\delta(x,y))  >-K|x-y|^2~~~\mbox{if $x_n=0
\,,$}\label{bcxnl}
\\
&& G(y,-D_y\psi_\delta(x,y)) <K|x-y|^2~~~\mbox{if $y_n=0
\,.$}\label{bcynl}
\end{eqnarray}
\end{prop}
   \medskip

{\noindent\bf Proof.} We only check \rec{bcxnl} the other case
being similar. \par

If $x\in\partial O$ then $d(x)=x_n=0$; moreover since $y\in
\overline O$, $d(y)=y_n\geq 0$. Thus we have
\begin{eqnarray*}
D_x\psi_\delta(x,y)&=& 2(x-y)-2\tilde
\varphi_\delta(\frac{x+y}{2},x-y)|x-y| e_n
\\
&+&  \left(D_X \tilde \varphi_\delta(\frac{x+y}{2},x-y)+D_Y \tilde
\varphi_\delta(\frac{x+y}{2},x-y)\right)  |x-y|y_n
\\
&-&  2 \tilde \varphi_\delta(\frac{x+y}{2},x-y)\frac{x-y}{|x-y|}y_n -
2A_1 y_n e_n-K_1\delta e_n.
\end{eqnarray*}
By taking in account that $n(x)=-e_n$   we have $$
D_x\psi_\delta(x,y)=p+q+r\,, $$ where
\begin{eqnarray*}
p&=&2(x-y)+2  \varphi (x,x-y)n(x)\, , \\q&=&2[\tilde
\varphi_\delta(\frac{x+y}{2},x-y)|x-y|- \varphi (x,x-y)]n(x)
\\&+&(D_X \tilde \varphi_\delta(\frac{x+y}{2},x-y)+D_Y \tilde
\varphi_\delta(\frac{x+y}{2},x-y))  |x-y|y_n\\
&-&2 \tilde \varphi_\delta(\frac{x+y}{2},x-y)\frac{x-y}{|x-y|}y_n \, ,\\
r&=&(2A_1 y_n +K_1\delta )n(x) \,.
\end{eqnarray*}
We first notice that, taking in account the definition of $\tilde
\varphi_\delta$, since $\varphi$ is homogeneous of degree $1$ with
respect to
$p$ and satisfies ${\bf (G2)}$ we have
\begin{eqnarray*}
\left| \tilde
\varphi_\delta(\frac{x+y}{2},x-y)-\varphi(x,\frac{x-y}{|x-y|})|x-y|
\right|&\le& K|x-y|(|x-y|+\Gamma_{\delta})\nonumber\\ &\le& K
[|x-y|^2+\delta+(y_n-x_n)]\nonumber\\ &\le &K  |x-y|^2 + K\delta +
K y_n \,. \nonumber
\end{eqnarray*}
Moreover, we have
$$
   |\tilde \varphi_\delta(\frac{x+y}{2},x-y)\frac{x-y}{|x-y|}y_n|\le
   K|x-y|^2.
   $$
By using the fact that $G(x,p)=0$ and  combining the above
estimates with the properties of $G$, we get
\begin{eqnarray*}
G(x,D_x\tilde\psi_\delta)&\ge& G(x,p+r)-K|q|\\ &\ge & \lambda(2A_1
y_n +K_1\delta)-K|x-y|^2-K\delta-Ky_n\,.
\end{eqnarray*}
Thus if $A_1$ and $K_1$ are large enough we get $$
G(x,D_x\tilde\psi_\delta)\ge -K|x-y|^2+K\delta\,. $$

   \hfill\sn

\medskip
\noindent{\bf Step 4.} Construction of the function  $\phi_0$ and its
main properties.\par

For $\bar C, A_2 >0$ large enough, we introduce the following function
\begin{equation}\label{phidelta}
\phi_0(x,y)=|x-x_0|^4-2\bar
C|x-x_0|^3(d(x)-d(x_0))+A_2(d(x)-d(x_0))^4\;.
\end{equation}

In order to check that the function $\Theta_\delta$ satisfy the
right boundary conditions we premise   the following Lemma whose
proof can be found in \cite{bdl03}. To formulate it, we use the
following notation~: for $p\in \R$ and $x\in \partial O $, $\tp :=
p- \langle p, n(x)\rangle n(x)$. $\tp $ represents the projection of $p$
on the tangent hyperplane to $\partial O$ at $x$.

\begin{lemma}\label{propG} Assume that {\bf (G1)} and {\bf (G2)} hold
and that, for some $x\in \partial O$,  and $\tilde p \in \R$, we
have $G(x,\tilde p) \leq 0$ (resp.  $G(x,\tilde p) \geq 0$), then
there exists a constant $\bar K $ (depending on $\nu$ and $K$ in
{\bf (G1)-(G2)}) such that, if $ \langle p, n(x)\rangle \leq -\bar K
\vert
\tp \vert$, then
$$ G(x,\tilde p+p) \leq 0\; ,$$ (resp. if $\langle p, n(x)\rangle \geq
\bar
K\vert \tp \vert$, then $$ G(x,\tilde p+p) \geq 0\; ).$$
\end{lemma}

The connection with $\phi_0$ is given by the following result.
\begin{lemma}\label{propphi}
For $\bar C$ large enough and for all $x\in \partial O$,  $\phi_0$
   satisfies
$$
\frac{\partial\phi_0}{\partial n}(x)\ge \bar K|{{\cal
T}(D\phi_0)}|\,,
$$
where   $\bar K>0$ is the constant appearing in Lemma \ref{propG}.
\end{lemma}
{\bf \noindent Proof.} We first compute the normal derivative of
$\phi_0$ and we use the usual property linking the distance function
and $n$; this yields
\begin{eqnarray*}
\langle D_x\phi_0 (x),n(x)\rangle &=& \langle
4|x-x_0|^2(x-x_0),n(x)\rangle+2\bar
C|x-x_0|^3\\ &=&4|x-x_0|^2(d(x_0)+O(|x-x_0|^2))+2\bar C|x-x_0|^3
\\ &\ge& 2\bar C|x-x_0|^3 +O(|x-x_0|^4)\,.
\end{eqnarray*}
On an other hand, we clearly have $ |{\cal T}(D\phi_0)| \leq
4|x-x_0|^3$.
Thus, if  $\bar K$ is the constant given in Lemma~\ref{propG},  by
choosing $\bar C$ large enough and $x$ close to $x_0$,  we have
\begin{eqnarray*}
   |{\cal T}(D\phi_0)|^{-1}\langle D_x\phi_0 (x),n(x)\rangle
   &\ge&
   \frac{2\bar C|x-x_0|^3 +O(|x-x_0|^4)}{4|x-x_0|^3}> \bar K\,.
\end{eqnarray*}
~~\hfill\sn

Now we can prove the
following result.
\begin{prop}\label{bc3}
If $|x-y|<<1$, then  we have
\begin{eqnarray}
&& G(x, D_x\Theta_{\delta}(x,y)) >0~~\mbox{if $x_n=0
\,,$}\label{hbcxn2}
\\
&& G(y,-D_y\Theta_{\delta}(x,y)) <0~~\mbox{if $y_n=0
\,.$}\label{hbcyn2}
\end{eqnarray}
   \end{prop}
{\noindent\bf Proof of Proposition \ref{bc3}.} Again we only check
the first property \rec{hbcxn2}.  First of all we observe that,
because of the assumption {\bf (G1)} and the property \rec{bcxnl}
of $\psi_\delta$  we have
\begin{eqnarray*}
G(x,D_x(e^{-\tilde K(y_n)} (\psi_\delta(x,y))^{\frac{\alpha}{2}})
&=&G(x,\tilde K e^{-\tilde Ky_n}
(\tilde\psi_\delta(x,y))^{\frac{\alpha}{2}}n(x)\\&+&\frac{
\alpha}{2} e^{-\tilde Ky_n}
(\tilde\psi_\delta(x,y))^{\frac{\alpha}{2}-1}D_x \psi_\delta(x,y))
\\&\ge&
\mu\tilde K e^{-\tilde K(y_n)} |x-y|^\alpha - K   e^{-\tilde
K(y_n)}|x-y|^\alpha\,.
\end{eqnarray*}
Thus, if $\tilde K$ is sufficiently large,  we obtain
$$G(x,D_x(e^{-\tilde K(x_n+y_n)}
(\psi_\delta(x,y))^{\frac{\alpha}{2}})>   K |x-y|^\alpha\; ,$$ for
some constant $  K$.
\par
Now by combining Lemma \ref{propG} and \ref{propphi} we get
\begin{eqnarray*}
G(x, D_x\Theta_{\delta}(x,y))&=&G(x,CD_x(e^{-\tilde K(x_n+y_n)}
(\psi_\delta(x,y))^{\frac{\alpha}{2}})+LD_x\phi_0)\\ &>&  K
|x-y|^\alpha. \end{eqnarray*}

And the result is proved.\hfill{$\sn$}

\par\bigskip \centerline{{\sc
Acknowledgements}}
\smallskip
    {The second author was
partially supported by M.I.U.R., project
   ``Viscosity, metric, and control theoretic methods for nonlinear
   partial differential equations'' and by G.N.A.M.P.A, project
``Equazioni alle derivate parzilai e teoria del controllo".}

\end{document}